\documentclass{article}

\usepackage{etex}
\usepackage{graphicx}
\usepackage{latexsym}
\usepackage{amsfonts}
\usepackage{amssymb}
\usepackage{amsmath}
\usepackage{verbatim}
\usepackage{url}
\usepackage[amsmath,standard, thmmarks]{ntheorem}
\usepackage{ntheorem}
\usepackage[margin=3cm]{geometry}
\usepackage{fancyhdr}

\newtheorem{thm}{Theorem}[section]
\newtheorem{cor}[thm]{Corollary}
\newtheorem{lem}[thm]{Lemma}
\newtheorem{prop}[thm]{Proposition}

\newtheorem{defn}[thm]{Definition}

\newcommand{\R}{\mathbb R}

\newcommand{\Z}{\mathbb Z}

\newcommand{\Sp}{\text{Sp\,}}

\newcommand{\To}{\longrightarrow}

\newcommand{\CC}{\mathcal{C}}

\newcommand{\E}{\mathcal{E}}
\newcommand{\z}{{\bf z}}

\newcommand{\D}{\mathcal{D}}
\newcommand{\Pent}{\mathcal{P}}
\newcommand{\V}{\mathcal{V}}

\newcommand{\hyp}{\mathbb H}

\DeclareMathOperator{\Par}{Par} \DeclareMathOperator{\Hyp}{Hyp}
\DeclareMathOperator{\Ell}{Ell}

 \DeclareMathOperator{\Isom}{Isom}
\DeclareMathOperator{\Axis}{Axis} \DeclareMathOperator{\Twist}{Tw}
\DeclareMathOperator{\Tr}{Tr} 
\DeclareMathOperator{\Fix}{Fix}

\DeclareMathOperator{\Top}{Top}

\numberwithin{equation}{section}

\begin{document}

\title{The hyperbolic meaning of the Milnor--Wood inequality}

\author{Daniel V. Mathews}%

\date{}

\maketitle
% ----------------------------------------------------------------

\begin{abstract}

We introduce a notion of the twist of an isometry of the hyperbolic plane. This twist function is defined on the universal covering group of orientation-preserving isometries of the hyperbolic plane, at each point in the plane. We relate this function to a function defined by Milnor and generalised by Wood. We deduce various properties of the twist function, and use it to give new proofs of several well-known results, including the Milnor--Wood inequality, using purely hyperbolic-geometric methods. Our methods express inequalities in Milnor's function as equalities, with the deficiency from equality given by an area in the hyperbolic plane. We find that the twist of certain products found in surface group presentations is equal to the area of certain hyperbolic polygons arising as their fundamental domains.

\end{abstract}

\tableofcontents

\section{Introduction}

\subsection{Overview}

In his 1957 paper \cite{Milnor}, Milnor introduced a function $\Theta: \widetilde{GL_2^+ \R} \To \R$ which is in a sense a ``rotation angle'' associated to elements of the universal covering group of the matrix group $GL_2^+ \R$. He proved that it satisfies the inequality
\[
\big| \Theta(\alpha \beta) - \Theta(\alpha) - \Theta(\beta) \big| < \frac{\pi}{2},
\]
i.e. is a \emph{quasimorphism}, and used it to prove a theorem regarding the existence of principal $GL_2^+ \R$ bundles over a closed oriented surface with a flat connection. This result was extended by Wood in \cite{Wood}, who defined a function $r: \widetilde{\Top^+ S^1} \to \R$, with similar properties; here $\Top^+ S^1$ is the group of orientation-preserving homeomorphisms of the circle and $\widetilde{\Top^+ S^1}$ its universal cover. Wood used this function $r$ to prove, \emph{inter alia}, a theorem regarding bundles over surfaces with structure group $\Top S^1$; in particular, when the structure group reduces to a totally disconnected subgroup.

One way to interpret the proofs of these theorems, broadly, is as follows. The function $\Theta$ or $r$ gives a measure of how far an element of $\tilde{G}$ (where $G$ is $GL_2^+ \R$ or $\Top^+ S^1$ or some other group) is from the origin. The quasimorphism property is used to show that a commutator of any two elements in $\tilde{G}$ cannot be ``too far'' from the origin. Since bundles over surfaces with flat connections (or totally disconnected structure group) are given by holonomy representations, understanding bundles of the desired type is essentially the same as understanding holonomy representations; and since an oriented surface has a standard presentation with one relator, namely a product of commutators, the understanding of commutators in $\tilde{G}$ gives results about the existence of such bundles.

The key result in these theorems, then, is what has become known as the \emph{Milnor--Wood inequality} (see e.g. \cite{Goldman88}), which expresses how far a product of commutators in $\tilde{G}$, which multiplies to $1 \in G$ (as required of a surface group representation) can stray from the identity. In particular, letting the lifts of $1 \in G$ to $\tilde{G}$ be $\{\z^m\}$, such a product of commutators is of the form $\z^m$; this $m$ is essentially the \emph{Euler class} of the representation and the content of the inequality is that this Euler class cannot exceed the Euler characteristic of the surface in magnitude.

The Milnor--Wood inequality is by now a classical result and has given rise to a vast array of applications and generalisations. For example: in the theory of Lorentz spacetimes of constant curvature \cite{Mess07}, circular groups \cite{Calegari04}, foliations \cite{Sullivan76, Eisenbud-Hirsch-Neumann}, contact geometry \cite{GiBundles}, and bounded cohomology \cite{Gromov82}.  It has been generalised to other Lie groups \cite{BG-PG_Maximal_surface_group_representations} and to general representations of lattices into Lie groups of Hermitian type \cite{BurgerIozzi}. Analogous results exist in higher-dimensional hyperbolic geometry \cite{Besson-Courtois-Gallot07} and for manifolds locally isometric to a product of hyperbolic planes \cite{Bucher-Gelander}. This is just a random sample and is by no means even an overview of the work which exists on the topic.

In this paper we present something far lower-powered, and restricted to Milnor's original case, but perhaps still of interest; we are surprised not to have found this idea in the existing literature. The present paper is concerned with $G = SL_2\R$; obviously $SL_2\R \subset GL_2^+ \R$ and $\widetilde{SL_2\R} = \widetilde{PSL_2\R} = \widetilde{\Isom^+\hyp^2}$. Milnor's $\Theta$ thus assins a number to a (lift to universal cover of a) hyperbolic isometry. We will give a hyperbolic-geometric interpretation of $\Theta$ by defining a function
\[
\Twist \; : \; \widetilde{PSL_2\R} \times \overline{\hyp}^2 \To \R,
\]
the ``twist angle'' of an $\tilde{\alpha} \in \widetilde{PSL_2\R}$ at a point $p \in \hyp^2$, which generalises $\Theta$. This function has interesting properties, including quasimorphism-type properties, which give a hyperbolic-geometric proof of the quasimorphism property of $\Theta$. Even better, we give an \emph{equality} in which the defect of $\Twist$ (and hence $\Theta$) from being a homomorphism is expressed as an area in the hyperbolic plane. Areas arise as deficiencies from additivity essentially because of the effect of negative curvature on parallel translation. Thus, we obtain a new proof of the Milnor--Wood inequality by pure hyperbolic-geometric methods.

We have several other applications. We use the function $\Twist$ to prove various relationships between surface group representations and areas in the hyperbolic plane. We interpret the twist of a commutator as the area of a hyperbolic pentagon, and indeed we can interpret the twist of any product occurring as a standard orientable surface group relator as an area of a polygon in the hyperbolic plane. We can also reprove some known results about hyperbolic isometries: which elements of $\widetilde{PSL_2\R}$ which can occur as commutators \cite{Wood, Eisenbud-Hirsch-Neumann, Goldman88}; relationships between types of commutators and their trace \cite{Goldman88, Goldman03}; and a cute result, as far as we know first appearing in \cite{Goldman03}, characterising isometries of hyperbolic type with intersecting axes in terms of the trace of their commutator. 

However, in our view the main new application of these methods is in a pair of subsequent papers, where we consider the question of which representations of the fundamental group of a surface are holonomy representations of hyperbolic structures, and cone-manifold structures of a certain type: see \cite{Me10MScPaper1, Me10MScPaper2}. Our methods here establish connections between the algebra of $PSL_2\R$ and hyperbolic geometry, which we use in those papers.

We finally note that $SL_2\R \cong \Sp(2)$, the group of $2 \times 2$ symplectic matrices, i.e. linear symplectomorphisms of $\R^2$ with the standard symplectic structure. Milnor's function $\Theta$ in this context is essentially the \emph{Maslov index} (see e.g. \cite[p. 48]{McDuffSalamon_Introduction}). We wonder if there are any further connections to symplectic geometry.

\subsection{Structure of this paper}

In section \ref{sec:twisting} we define the notion of twist. This first requires some preliminaries on $\widetilde{PSL_2\R}$, which occupy sections \ref{subsec:PSL_and_PSL} to \ref{subsec:traces}. In section \ref{sec:properties_of_twisting} we establish various properties of our twist function. In section \ref{sec:Milnors_function} we recall the definition of Milnor's function $\Theta$, we relate it to twisting, and deduce various properties.

\subsection{Acknowledgments}

This paper forms one of several papers arising from the author's Masters thesis \cite{Mathews05}, completed at the University of Melbourne under Craig Hodgson, whose advice and suggestions have been highly valuable. It was completed during the author's postdoctoral fellowship at the Universit\'{e} de Nantes, supported by a grant ``Floer Power'' from the ANR.

\section{Twisting in the hyperbolic plane}
\label{sec:twisting}

Everything in sections \ref{subsec:PSL_and_PSL} to \ref{subsec:traces} has been known for a long time: see, e.g. \cite{Goldman88}. Although the idea is very basic, it appears that the notion of the twist function which we define in section \ref{sec:deriv_isoms} is new.

\subsection{$PSL_2\R$ and $\widetilde{PSL_2\R}$}
\label{subsec:PSL_and_PSL}

Fix a basepoint $y_0$ in $\hyp^2$ and unit tangent vector $u_0 \in UT_{y_0} \hyp^2$. An orientation-preserving hyperbolic isometry is uniquely determined by the image of $u_0$, and we may identify the unit tangent bundle $UT\hyp^2$ with the orientation-preserving isometry group $PSL_2\R$. Topologically $PSL_2\R \cong \R^2 \times S^1$; let $p_1$ be the projection map $PSL_2\R \To \hyp^2$. 

Let $p_2: \widetilde{PSL_2\R} \To PSL_2\R$ be the universal cover of $PSL_2\R$; see \cite{Goldman_thesis, Goldman88} for further details. Clearly $\pi_1(PSL_2\R) \cong \Z$. An element $\tilde{x} \in \widetilde{PSL_2\R}$ is hyperbolic, elliptic or parabolic accordingly as is $p_2(\tilde{x}) \in PSL_2\R$. We can consider $\widetilde{PSL_2\R}$ as $\hyp^2$, with an $\R$ fibre above each point, covering the circle of unit tangent vectors.
\[
   \widetilde{PSL_2\R} \stackrel{p_2}{\To} PSL_2\R \cong UT\hyp^2 \stackrel{p_1}{\To} \hyp^2.
\]

We can also consider elements of $\widetilde{PSL_2\R}$ as homotopy classes of paths in $UT\hyp^2$ starting at the basepoint. Since the basepoint is arbitrary,
\emph{every} path $c: [0,1] \To UT\hyp^2$ determines a unique element of $\widetilde{PSL_2\R}$, which we also denote $c$, abusing notation. The projection of $c$ to $PSL_2\R$ is the orientation-preserving isometry sending $c(0)$ to $c(1)$. An $\alpha \in PSL_2\R$ has countably infinitely many lifts to $\widetilde{PSL_2\R}$. These all represent paths in $UT\hyp^2$ between the same start and end tangent vectors. However these paths will differ according to the number of times that the tangent vectors spin as the path is traversed. The lifts of the identity $1 \in PSL_2\R$ form an infinite cyclic group $\{\z^n :  n \in \Z \}$, where $\z$ is the homotopy class of the path $c(t) = (y_0, e^{2\pi int} u_0 )$. Note $\z$ commutes with every element of $\widetilde{PSL_2\R}$; in fact $\z$ generates the centre of $\widetilde{PSL_2\R}$.

\subsection{Regions in $\widetilde{PSL_2\R}$}

While every element has infinitely many lifts, some lifts are simpler than others. For instance, the identity in $\widetilde{PSL_2\R}$ is the ``simplest'' lift of the identity in $PSL_2\R$.

If $\alpha \in PSL_2\R$ is hyperbolic then it translates by distance $d_\alpha$ along $\Axis \alpha$. Let $c(t) \in PSL_2\R$ be the translation of (signed) hyperbolic distance $td_\alpha$ along $\Axis \alpha$; then $c: \R \To PSL_2\R$ is a homomorphism with $c(1)=\alpha$, in fact the only homomorphism with this property. The path $c|_{[0,1]}$ in $PSL_2\R$ gives an element $\tilde{\alpha}$ of $\widetilde{PSL_2\R}$ which we take as our preferred or simplest lift. This lift can be considered a path of unit tangent vectors, which travels along $\Axis \alpha$ at speed $d$, always pointing along $\Axis \alpha$ in the direction of translation.

\begin{figure}[tbh]
\centering
\includegraphics[scale=0.4]{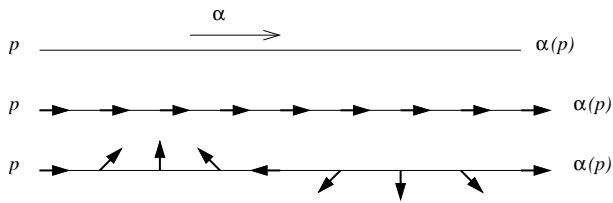}
\caption{An isometry $\alpha \in PSL_2\R$; the simplest lift of
$\alpha$; a different lift of $\alpha$.}
\end{figure}

Similar considerations apply to parabolic isometries. A parabolic $\alpha \in PSL_2\R$ translates along some horocycle $h_\alpha$ (not unique); endowing $h_\alpha$ with a Euclidean metric, let $\alpha$ translate by Euclidean distance $d$. Letting $c(t) \in PSL_2\R$ be the parabolic translating $td$ along $h_\alpha$ then $c: \R \To PSL_2\R$ is the unique homomorphism with $c(1)=\alpha$, and $c|_{[0,1]}$ gives a preferred lift $\tilde{\alpha} \in \widetilde{PSL_2\R}$. This $\tilde{\alpha}$ can be considered a path of tangent vectors travelling along and pointing along $h_\alpha$ at speed $d$ for time $1$.

However the situation for elliptic $\alpha \in PSL_2\R$ is different: there are infinitely many homomorphisms $c: \R \To PSL_2\R$ with $c(1)=\alpha$. Let $\alpha$ rotate by angle $\theta$ (mod $2\pi$). Then the lifts of $\alpha$ are rotations by angles $\theta + 2\pi \Z$. From this viewpoint there are two simplest lifts of $\alpha$, those with rotation angle lying in $(0, 2\pi)$ and $(-2\pi,0)$: a simplest anticlockwise and clockwise lift.

Denote the sets of simplest lifts of hyperbolics and parabolics $\Hyp_0$ and $\Par_0$ respectively. Let $\Hyp_n = \z^n \Hyp_0$ and $\Par_n = \z^n \Par_0$, so the hyperbolic (resp. parabolic) elements of $\widetilde{PSL_2\R}$ are $\sqcup_n \Hyp_n$ (resp. $\sqcup_n \Par_n$). We may consider $\tilde{\alpha} \in  \Hyp_n$ as a translation along $\Axis \alpha$ with an added twist of $2n\pi$. We may further distinguish between $\Par_n^+$ and $\Par_n^-$, the rotations about points at infinity whose projections to $PSL_2\R$ are anticlockwise and clockwise respectively.

As for elliptics, let the set of simplest anticlockwise and clockwise lifts be $\Ell_1$ and $\Ell_{-1}$ respectively. For $n>0$ let $\Ell_n = \z^{n-1} \Ell_1$ and $\Ell_{-n} = \z^{-n+1} \Ell_{-1}$. (So $\Ell_0$ is not defined and $\z \Ell_{-1} = \Ell_1$.) For $n>0$ (resp. $n<0$), $\Ell_n$ consists of all rotations through angles between $2\pi(n-1)$ and $2\pi n$ (resp. between $2\pi n$ and $2\pi (n+1)$).

\begin{figure}[tbh]
\centering
\includegraphics[scale=0.6]{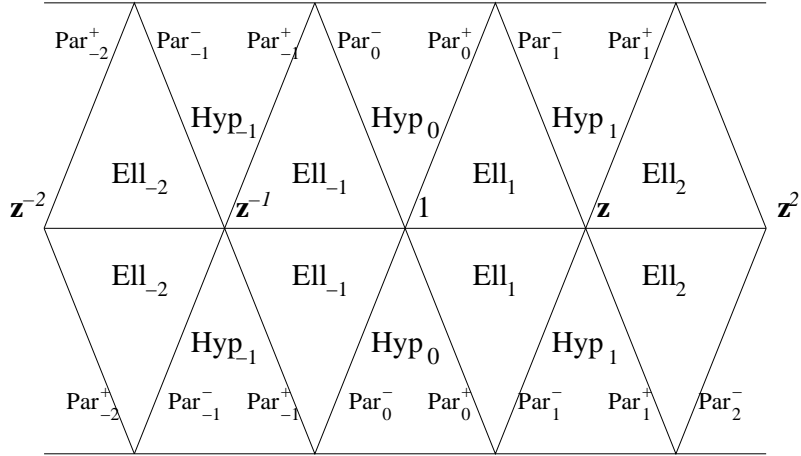}
\caption[Schematic diagram of $\widetilde{PSL_2\R}$]{Schematic
diagram of $\widetilde{PSL_2\R}$.}
\label{universalcover}
\end{figure}

Considering that hyperbolics and elliptics form 3-dimensional subspaces of the 3-dimensional $PSL_2\R$, with common 2-dimensional boundary the space
of parabolic elements, we may draw a schematic diagram of $\widetilde{PSL_2\R}$ as in figure \ref{universalcover}.

\begin{lem}
\label{commutator_lift}
    Let $\alpha, \beta \in PSL_2\R$. Then $[\alpha,\beta]$
    has a
    well-defined lift to $\widetilde{PSL_2\R}$. That is, any two
    sets of lifts $\tilde{\alpha}_1, \tilde{\beta}_1$ and $\tilde{\alpha}_2,
    \tilde{\beta}_2$ satisfy $[\tilde{\alpha}_1, \tilde{\beta}_1] =
    [\tilde{\alpha}_2, \tilde{\beta}_2]$.
\end{lem}

\begin{Proof}
    Let $\tilde{\alpha}_2 = \z^a \tilde{\alpha}_1$, $\tilde{\beta}_2 = \z^b
    \tilde{\beta}_1$. Since $\z$ is central,
    \begin{align*}
        [\tilde{\alpha_2}, \tilde{\beta_2}] &= \tilde{\alpha_2} \tilde{\beta_2}
        \tilde{\alpha_2}^{-1} \tilde{\beta_2}^{-1}
        = \z^a \tilde{\alpha_1} \z^b
        \tilde{\beta_1} \tilde{\alpha_1}^{-1} \z^{-a} \tilde{\beta_1}^{-1}
        \z^{-b}% \\
        = \tilde{\alpha_1} \tilde{\beta_1} \tilde{\alpha_1}^{-1}
        \tilde{\beta_1}^{-1}
        = [\tilde{\alpha_1}, \tilde{\beta_1}].
    \end{align*}\
\end{Proof}

\subsection{Traces in $\widetilde{PSL_2\R}$}
\label{subsec:traces}

As $\widetilde{PSL_2\R}$ covers $SL_2\R$, there is a well-defined trace on $\widetilde{PSL_2\R}$. For all elliptic regions, the trace lies in $(-2,2)$; in the various other regions of $\widetilde{PSL_2\R}$ it takes values as follows.

\begin{lem}[Trace lemma]
\label{trace_regions}
\[
        \Tr \left( \z^n \right) = (-1)^n \cdot 2, \quad
        \Tr \left( \Par_n \right) = (-1)^n \cdot 2, \quad
        \Tr \left( \Hyp_n \right) = \left\{ \begin{array}{ll} (2,
        \infty) & \text{$n$ even} \\ (-\infty, -2) & \text{$n$ odd.}
        \end{array} \right.
\]
\end{lem}

\begin{Proof}
Consider the matrix
\[
  E(\theta) = \begin{pmatrix} \cos \theta & \sin \theta \\  - \sin \theta & \cos \theta \end{pmatrix}.
\]
Now in the upper half plane $E(\theta)$ is elliptic, a rotation of $2\theta$ about $i$. Thus $E(n \pi) = \z^n$, and hence $\z^n$ projects to $(-1)^n \in SL_2\R$. From this the first claim follows. The trace of an element of $\widetilde{PSL_2\R}$ is $\pm 2$ if and only if it is a power of $\z$, or parabolic. Now $\Tr$ is a continuous function and, considering the topology of $\widetilde{PSL_2\R}$, $\Par_n \cup \{\z^n\}$ is connected. Thus $\Tr (\Par_n) = (-1)^n$. As $\Hyp_n$ is connected and bounded by $\Par_n$ and $\z^n$, on which $\Tr = (-1)^n \cdot 2$, the final claim follows.
\end{Proof}

\subsection{Definition of twist}

\label{sec:deriv_isoms}

Define the \emph{twist of a vector field along a curve} as follows. Consider a smooth curve $c: [0,1] \To \hyp^2$ and a smooth unit tangent vector field $\V: [0,1] \To UT\hyp^2$, $p_1 \circ \V = c$ (recall $p_1$ is the projection $UT\hyp^2 \To \hyp^2$). Consider the velocity vector field $\frac{dc}{dt}$ along $c$, which we may rescale to a unit vector field $\hat{c}: [0,1] \To UT\hyp^2$. Consider the angle $\theta(t)$ (measured anticlockwise) from $\hat{c}(t)$ to $\V(t)$. We have many choices for $\theta(0)$ (differing by $2\pi \Z$), but choosing $\theta(0)$ arbitrarily determines continuous $\theta$ completely; $\theta(1) - \theta(0)$ is independent of this choice, and is the twist of $\V$ along $c$.

\begin{figure}[tbh]
\centering
\includegraphics[scale=0.4]{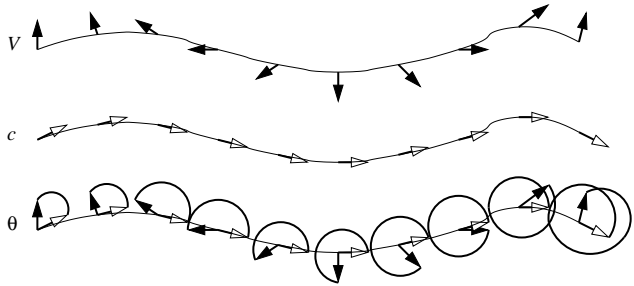}
\caption{The twist of a vector field along a curve.}
\end{figure}

Now given $y \in \hyp^2$ and $\tilde{\alpha} \in \widetilde{PSL_2\R}$ we define the \emph{twist of $\tilde{\alpha}$ at $y$}, denoted $\Twist(\tilde{\alpha},y)$. Let $\tilde{\alpha}$ project to $\alpha \in PSL_2\R$. Let $c: [0,1] \To \hyp^2$ be a constant speed (possibly $0$) geodesic from  $y$ to $\alpha(y)$. There is a vector field $\V: [0,1] \To UT\hyp^2$ along $c$ which lies in the homotopy class of $\tilde{\alpha}$. Then $\Twist(\tilde{\alpha},y)$ is the twist of $\V$ along $c$.

That is, $\Twist(\tilde{\alpha},y)$ describes how the tangent vector at $y$ is moved by $\tilde{\alpha}$, compared to parallel translation along the geodesic from $y$ to $\tilde{\alpha}(y)$. It is clear this does not depend on the choice of $\V$. For $\alpha \in PSL_2\R$, define $\Twist(\alpha, y)$ the same way, except angles are taken modulo $2\pi$.

As an aside, we note that this method of obtaining a rotation angle from an element of $\widetilde{PSL_2\R}$ is not so different from Wood's method in \cite{Wood}. Wood regards an element of $SL_2\R$ as acting on $\R^2$ by a linear transformation, and hence on the $S^1$ of oriented lines through the origin. Thus there is an inclusion $SL_2\R \subset \Top^+ S^1$. This is equivalent to the action of a hyperbolic isometry on the circle at infinity. But here we regard the isometry as acting on the $S^1$ of unit tangent vectors at a given point; such unit tangent vectors of course correspond bijectively with the circle at infinity, but different points give different bijections. Our twist is the action of an isometry on unit tangent $S^1$'s, where the $S^1$'s at a point and its image are related by parallel translation. Wood's function on $\widetilde{\Top^+ S^1}$ involves an integral and hence the measure on the circle at infinity; an isometry however alters this measure.

\section{Properties of twisting}
\label{sec:properties_of_twisting}

\subsection{Types of isometries}

One can easily verify the following properties of the twist.
\begin{itemize}
\item 
For hyperbolic $\alpha \in PSL_2\R$ and $y \in \Axis \alpha$, $\Twist(\alpha,y) = 0$ (mod $2\pi$). For $\tilde{\alpha} \in \Hyp_0$, $\Twist(\tilde{\alpha},y) \in (-\pi,\pi)$. The twist is constant along curves of constant distance $h$ from $\Axis \alpha$. For each $\theta \in (-\pi,\pi)$ there is precisely one $h$ for which the curve at distance $h$ is the locus of points $y$ with $\Twist(\tilde{\alpha},y) = \theta$.

\item
For $\tilde{\alpha} \in \Par_0$, $\Twist(\tilde{\alpha},y)$ is constant along horocycles about $\Fix \alpha$. If $\tilde{\alpha} \in \Par_0^+$ (resp. $\Par_0^-$) then $\Twist(\tilde{\alpha},y) \in (0,\pi)$ (resp. $(-\pi,0)$). For horocycles close to $\Fix \alpha$, the twist is close to $0$. For each $\theta \in (0, \pi)$ (resp. $(-\pi, 0)$) there is precisely one horocycle which is the locus of points $y$ with $\Twist(\tilde{\alpha},y) = \theta$.

\item
For elliptic $\tilde{\alpha}$, $\Twist(\tilde{\alpha},y)$ is constant along hyperbolic circles centred at $\Fix \alpha$. Take $\tilde{\alpha} \in \Ell_1$ for convenience, so $\tilde{\alpha}$ rotates by angle $\psi \in (0, 2\pi)$. So $\Twist(\tilde{\alpha}, \Fix \alpha) = \psi$. If $\psi \in (0,\pi)$ then $\Twist(\tilde{\alpha},y)$ always lies in $[\psi, \pi)$; for each $\theta \in (\psi, \pi)$ there is precisely one hyperbolic circle centred at $\Fix \alpha$ with is the locus of $y$ with $\Twist(\tilde{\alpha},y) = \theta$. If $\psi = \pi$ then $\alpha$ is a half turn and $\Twist(\tilde{\alpha}, y) = \pi$ for all $y$. If $\psi \in (\pi, 2\pi)$ then $\Twist(\tilde{\alpha},y)$ always lies in $(\pi,\psi]$ and for each $\theta \in (\psi, \pi)$ there is precisely one hyperbolic circle centred at $\Fix \alpha$ which is the locus of $y$ with $\Twist(\tilde{\alpha},y) = \theta$.
\end{itemize}

The values of the twist for all $\Hyp_n$, $\Par_n$ and $\Ell_n$ follow obviously from the above.
\begin{prop}\
\label{twist_bounds}
    \begin{align*}
        \Twist(\Hyp_n, \hyp^2) &= \Big( \left( 2n-1 \right) \pi,
        \left( 2n + 1 \right) \pi \Big) \\
        \Twist(\Par_n, \hyp^2) &= \Big( \left( 2n-1 \right) \pi,
        \left( 2n + 1 \right) \pi \Big) \\
        \Twist(\Ell_n, \hyp^2) &= \left\{ \begin{array}{ll} \Bigl( (2n-2)\pi, 2n\pi \Bigr) & \text{
        for $n>0$} \\
        \Bigl( -2|n|\pi, (-2|n|+1)\pi \Bigr) &
        \text{ for $n<0$} \end{array} \right.
    \end{align*}\
    \qed
\end{prop}

\subsection{Extension to infinity}

We have defined $\Twist \left( \tilde{\alpha}, p \right)$ for $p \in \hyp^2$. We can extend the definition to $p \in \overline{\hyp}^2$, with the circle at infinity $S_\infty^1$ adjoined. However we pay the price that, at least if $S_\infty^1$ is endowed with the usual topology of $S^1$, then $\Twist$ is not continuous on $S_\infty^1$.

Extending the definition is simple enough. Take $\tilde{\alpha} \in \widetilde{PSL_2\R}$ and $p \in S_\infty^1$, We note that for any geodesic $l$ with an endpoint at $p$, if we take points $x \in l$ approaching $p$ then $\Twist \left( \tilde{\alpha}, x \right)$ approaches a limit; and the limit is independent of choice of $l$. In particular:
\begin{itemize}
\item
For $\tilde{\alpha} \in \Hyp_0$, $\Twist(\tilde{\alpha},p) \in \{ - \pi, 0, \pi \}$. The twist is $0$ at the two endpoints of $\Axis \alpha$ and $\pm \pi$ on the two open arcs of $S_1^\infty$ to either side. 
\item
For $\tilde{\alpha} \in \Par_0^+$ (resp. $\Par_0^-$), $\Twist(\tilde{\alpha},p)$ is $0$ when $p = \Fix \alpha$ and $\pi$ (resp. $-\pi$) otherwise. 
\item
For $\tilde{\alpha} \in \Ell_1$ (resp. $\Ell_{-1}$) then $\Twist(\tilde{\alpha},p) = \pi$ (resp. $-\pi$) for all $p \in S_\infty^1$. 
\end{itemize}
For $\tilde{\alpha}$ in $\Hyp_n$, $\Par_n^\pm$ or $\Ell_n$ in general, we adjust by the appropriate multiple of $2\pi$. In particular, on $S_\infty^1$, the twist is always an integer multiple of $\pi$.

\subsection{Parallel translation and curvature}

Recall that curvature is the effect on tangent vectors induced by parallel translation around a loop. As the hyperbolic plane has constant curvature $-1$, parallel translation around a loop gives a rotation on a tangent vector equal to the negative area enclosed. 

\begin{figure}[tbh]
\centering
\includegraphics[scale=0.5]{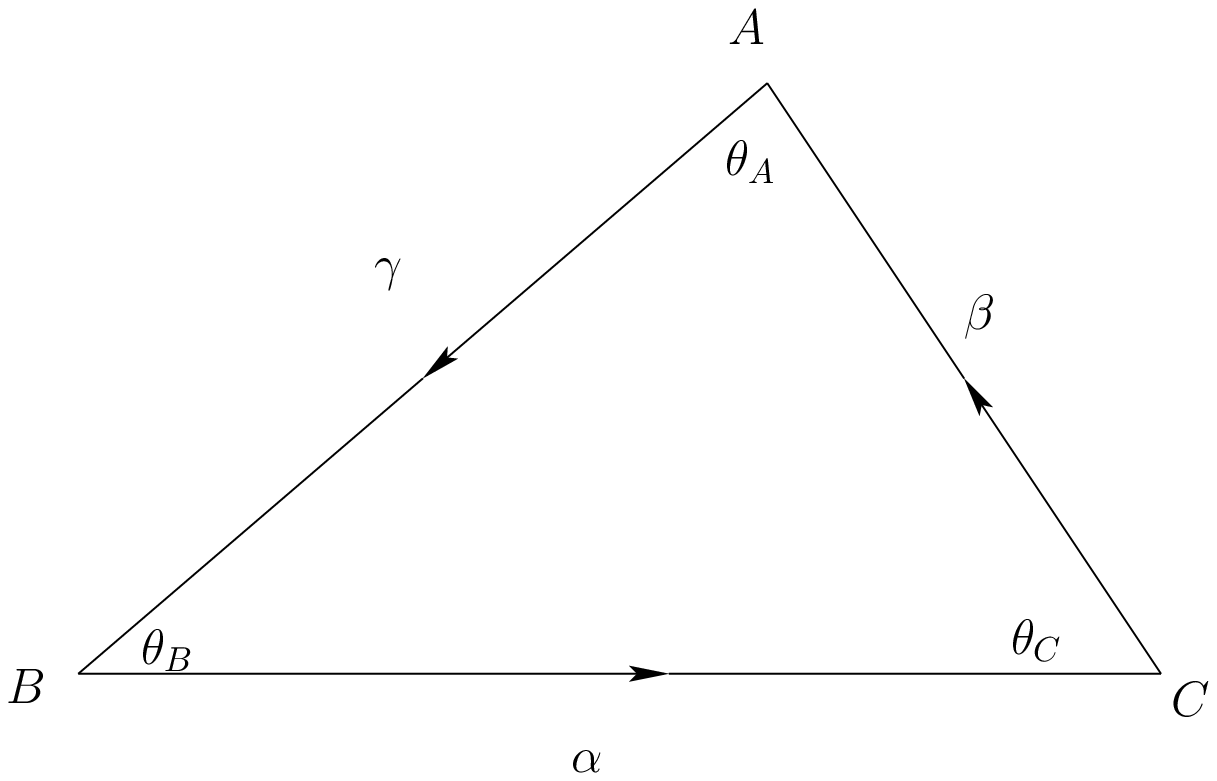}
\caption{Translation around a hyperbolic triangle.}
\label{fig:a}
\end{figure}

Let $ABC$ be a hyperbolic triangle shown in  figure \ref{fig:a}; let $\alpha, \beta, \gamma \in \Isom^+ \hyp^2$ respectively be the hyperbolic translations along axes $BC,CA,AB$ and translating $B \mapsto C \mapsto A \mapsto B$. Let $\tilde{\alpha}, \tilde{\beta}, \tilde{\gamma} \in \Hyp_0$ be their simplest lifts. Then the composition $\tilde{\gamma} \circ \tilde{\beta} \circ \tilde{\alpha}$ is parallel translation around $ABC$, hence a rotation of signed angle $\theta_A + \theta_B + \theta_C - \pi = - \Delta$, where $\Delta$ is the area of $ABC$. Hence, parallel translation from $B$ to $A$ is equivalent to first rotating an angle of $\Delta$ at $B$, then parallel translating $B \mapsto C \mapsto A$.

Let now $\tilde{\alpha}, \tilde{\beta}$ be \emph{any} elements of $\widetilde{PSL_2\R}$ covering hyperbolic isometries of $\hyp^2$ which take $B \mapsto C$ and $C \mapsto A$ respectively. Then $\tilde{\beta} \tilde{\alpha}$ takes $B \mapsto A$ and $\Twist(\tilde{\beta} \tilde{\alpha}, B)$ is given by the twist of $\tilde{\beta} \tilde{\alpha}$ along $BA$ relative to parallel translation. If we instead measured the twist of $\tilde{\beta} \tilde{\alpha}$ relative to parallel translation along $B \mapsto C \mapsto A$, i.e. $\Twist(\tilde{\alpha},B) + \Twist(\tilde{\beta},C)$, the answer must differ by $\Delta$. Hence we have the following result. Here and below we write $\Delta[A,B,C]$ to denote the signed area of the triangle with ordered vertices $A,B,C$, and use $\Delta$ in general to signify area. Taking a limit of points going to infinity, the result also holds for ideal points.
\begin{lem}[Composition lemma]
\label{lem:composition}
For any $\tilde{\alpha}, \tilde{\beta} \in \widetilde{PSL_2\R}$ and any $p \in \overline{\hyp}^2$,
\[
\Twist \left( \tilde{\beta} \tilde{\alpha}, p \right) = \Twist \left( \tilde{\alpha}, p \right) + \Twist \left( \tilde{\beta}, \alpha p \right) - \Delta[p,\alpha p, \beta \alpha p],
\]
where $\alpha,\beta \in \Isom^+ \hyp^2$ are the images of $\tilde{\alpha},\tilde{\beta}$.
\qed
\end{lem}

Thus, the failure of $\Twist$ to be linear is a manifestation of negative curvature, and the defect is the area of the triangle around which vectors are translated. The defect is clearly bounded as hyperbolic triangles have area less than $\pi$. For $p \in \hyp^2$ then,
\[
\left| \Twist \left( \tilde{\beta} \tilde{\alpha}, p \right) - \Twist \left( \tilde{\alpha}, p \right) - \Twist \left( \tilde{\beta}, \alpha p \right) \right| < \pi.
\]
For $p \in S_\infty^1$ the inequality holds but is not strict.

\subsection{Conjugation and addition}

\begin{lem}[Conjugation lemma]
\label{lem:conjugation}
For any $\tilde{\alpha}, \tilde{\beta} \in \widetilde{PSL_2\R}$ and $p \in \overline{\hyp}^2$, $\Twist (\tilde{\alpha},p) = \Twist ( \tilde{\beta} \tilde{\alpha} \tilde{\beta}^{-1}, \beta p)$.
\end{lem}

\begin{Proof}
Consider $\Twist (\tilde{\alpha},p)$, i.e. $\tilde{\alpha}$ as a path of unit tangent vectors along the geodesic from $p$ to $\alpha p$, compared to parallel translation. Now translate the whole situation by the isometry $\beta$.
\end{Proof}
When $\tilde{\beta} = \tilde{\alpha}^n$, for $n \in \Z$, this becomes
\begin{equation}
\Twist \left( \tilde{\alpha}, p \right) = \Twist \left( \tilde{\alpha}, \alpha^n p \right), 
\label{eq:translation} 
\end{equation}
which is clear geometrically: $\alpha^n p$ lies on the same constant distance curve from $\Fix \alpha$ or $\Axis \alpha$ as $p$.

\begin{lem}[Inverse lemma]
\label{lem:inverse}
For any $\tilde{\alpha} \in \widetilde{PSL_2\R}$ and $p \in \overline{\hyp}^2$,
\[
\Twist \left( \tilde{\alpha}, p \right) = - \Twist \left( \tilde{\alpha}^{-1}, p \right)
\]
\end{lem}

\begin{Proof}
Reversing the path of unit tangent vectors of $\tilde{\alpha}$ gives immediately $\Twist(\tilde{\alpha},p) = -\Twist(\tilde{\alpha}^{-1}, \alpha p)$. Now apply the previous corollary.
\end{Proof}
Apply composition lemma \ref{lem:composition} to the product $\tilde{\beta} \tilde{\alpha} = \left( \tilde{\beta} \tilde{\alpha} \tilde{\beta}^{-1} \right) \tilde{\beta}$
\[
\Twist \left( \tilde{\beta} \tilde{\alpha}, p \right) = \Twist \left( \tilde{\beta}, p \right) + \Twist \left( \tilde{\beta} \tilde{\alpha} \tilde{\beta}^{-1}, \beta p \right) - \Delta[p, \beta p, \beta \alpha p],
\]
then apply conjugation lemma \ref{lem:conjugation}. This gives a result like \ref{lem:composition}, but now all based at the same $p \in \hyp^2$.

\begin{lem}[Addition lemma]
\label{lem:addition}
For all $\tilde{\alpha}, \tilde{\beta} \in \widetilde{PSL_2\R}$ and $p \in \overline{\hyp}^2$,
\[
\Twist \left( \tilde{\beta} \tilde{\alpha}, p \right) = \Twist \left( \tilde{\beta}, p \right) + \Twist \left( \tilde{\alpha}, p \right) - \Delta[p, \beta p, \beta \alpha p].
\]
\qed
\end{lem}

This ``addition formula'' for $\Twist$ describes the twist of a product in terms of the twist of the factors, all at the same point. We immediately obtain a quasimorphism property of $\Twist$: for $p \in \hyp^2$,
\begin{equation}
\left| \Twist \left( \tilde{\beta} \tilde{\alpha}, p \right) - \Twist \left( \tilde{\beta}, p \right) - \Twist \left( \tilde{\alpha}, p \right) \right| < \pi.
\end{equation}
If $p \in S_1^\infty$ the inequality is not strict.

\subsection{Pentagons and commutators}

\begin{defn}
    Let $\alpha, \beta \in PSL_2\R$ and $p \in \hyp^2$. Then the geodesic
    pentagon in $\hyp^2$ obtained by joining the segments
    \[
	p \To \alpha^{-1} \beta^{-1} \alpha \beta p \To \beta p \To \alpha \beta p \To \beta^{-1} \alpha \beta p \To p
    \]
    is called \emph{the pentagon generated by $\alpha, \beta$ at $p$} and is
    denoted $\Pent(\alpha, \beta; p)$.
\end{defn}

Note that $\Pent(\alpha, \beta; p)$ may intersect itself; its vertices may coincide; it need not even bound an immersed disc. It is simply five geodesic line segments, possibly degenerate, in $\hyp^2$; if these all have nonzero length we say $\Pent(\alpha,\beta;p)$ is \emph{nondegenerate}. Denote the vertices
\[
    p_0 = p, \quad p_1 = \beta p, \quad p_2 = \alpha \beta p, \quad p_3 = \beta^{-1} \alpha \beta p, \quad p_4 = \alpha^{-1} \beta^{-1} \alpha \beta p.
\]
If $\Pent(\alpha,\beta;p)$ is nondegenerate and bounds an immersed disc, denote the interior angles of the pentagon $\theta_0, \ldots, \theta_4$. We may also denote a polygon as the sequence of vertices; we write $\Pent(\alpha,\beta;p) = [p_0, p_4, p_1, p_2, p_3]$. 

\begin{figure}[tbh]
\centering
\includegraphics[scale=0.5]{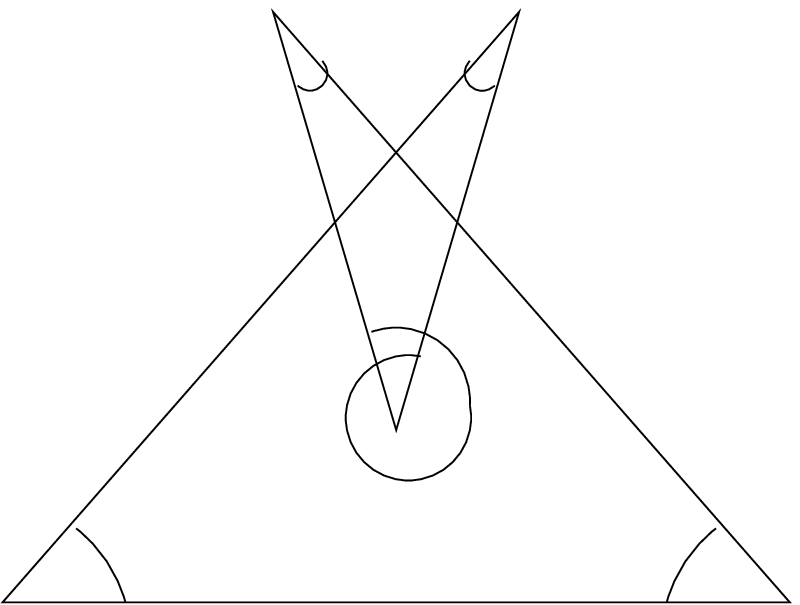}
\caption{$\Pent$ may bound an immersed but not embedded disc.}
\label{fig:26}
\end{figure}

If $\Pent(\alpha,\beta;p)$ bounds an embedded disc, then it has a well-defined area $\Delta[\Pent(\alpha,\beta;p)]$. This area is signed according to the boundary orientation $p_0 \rightarrow p_4 \rightarrow p_1 \rightarrow p_2 \rightarrow p_3 \rightarrow p_0$. The same can be done even if $\Pent(\alpha, \beta; p)$ bounds an immersed disc (see figure \ref{fig:26}); for instance by cutting into smaller embedded discs.

\begin{prop}[Commutator pentagon]
\label{prop:pentagon_twist}
If $\Pent(\alpha,\beta;p)$ is nondegenerate and bounds an immersed disc, then
\[
  \Twist \left( [\tilde{\alpha}^{-1},\tilde{\beta}^{-1}],p \right) = \Delta[\Pent(\alpha,\beta;p)].
\]
\end{prop}
(Recall by lemma \ref{commutator_lift}, the commutator is independent of choice of lifts $\tilde{\alpha}, \tilde{\beta}$.)

\begin{Proof}
Without loss of generality assume $\Delta > 0$. Consider figure \ref{fig:52}. Consider a unit tangent vector $(p,u)$ at $p$ pointing along the geodesic to $p_4$. Follow (``chase'') the image of this vector under $D\beta$, $D\alpha$, $D\beta^{-1}$ and $D\alpha^{-1}$ to obtain unit tangent vectors $(p_i, u_i)$ at
each $p_i$. Note that $\alpha$ takes the segment $p_4 \rightarrow p_1$ to the segment $p_3 \rightarrow p_2$ and $\beta$ takes $p_0 \rightarrow p_3$ to $p_1 \rightarrow p_2$; we use these two facts repeatedly.

\begin{figure}[tbh]
\centering
\includegraphics[scale=0.5]{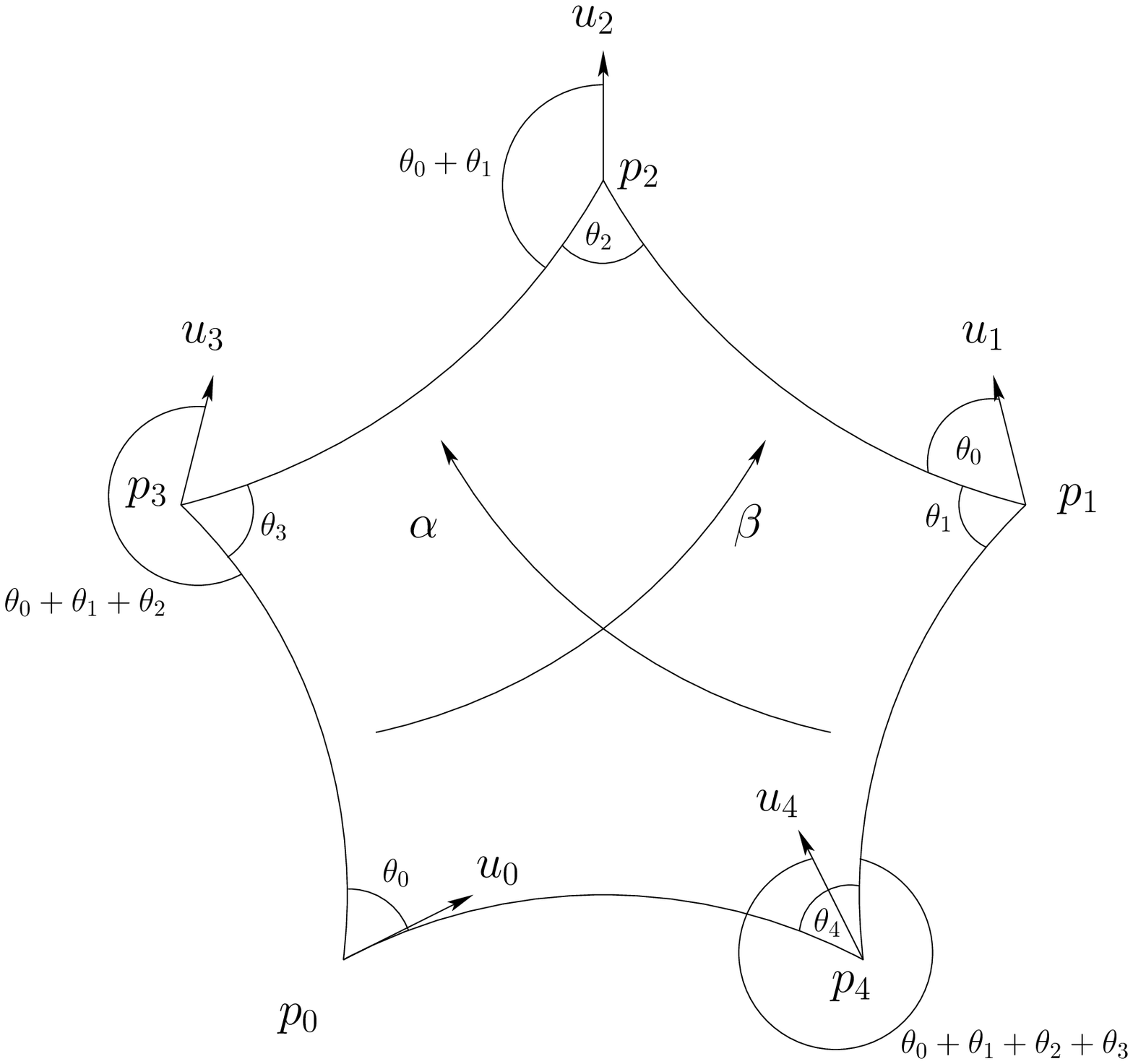}
\caption{A unit vector chase.} \label{fig:52}
\end{figure}

Now $(p_0,u_0) = (p,u)$ is based at $p$ and points $\theta_0$ clockwise of $p_0 \rightarrow p_3$; hence $(p_1, u_1) = D\beta (p,u)$ is based at $p_1$ and points $\theta_0$ clockwise of $p_1 \rightarrow p_2$. But $(p_1, u_1)$ points $\theta_0 + \theta_1$ clockwise of $p_1 \rightarrow p_4$, hence $(p_2, u_2) = D\alpha (p_1, u_1)$ points $\theta_0 + \theta_1$ clockwise of $p_2 \rightarrow p_3$. Then $(p_2, u_2)$ lies $\theta_0 + \theta_1 + \theta_2$ clockwise of $p_2 \rightarrow p_1$, hence $(p_3, u_3) = D\beta^{-1}(p_2,u_2)$ lies the same angle clockwise of $p_3 \rightarrow p_0$. Finally, $(p_3,u_3)$ lies $\theta_0 + \theta_1 + \theta_2 + \theta_3$ clockwise of $p_3 \rightarrow p_2$, and hence $(p_4,u_4) = D[\alpha^{-1},\beta^{-1}](p,u)$ lies the same angle clockwise of $p_4 \rightarrow p_1$. Thus $D[\alpha^{-1},\beta^{-1}](p,u)$ lies $\theta_0 + \theta_1 + \theta_2 + \theta_3 + \theta_4$ clockwise of $p_4 \rightarrow p_0$. 

This immediately shows that $\Twist ([\alpha^{-1},\beta^{-1}],p)= 3 \pi - \sum \theta_i = \Delta[\Pent(\alpha,\beta;p)]$, modulo $2\pi$. Choosing particular lifts $\tilde{\alpha}, \tilde{\beta}$ of $\alpha, \beta$ we may see that this is true over the real numbers. For this we use the following lemma, which is straightforward, although in the elliptic case perhaps the reader might draw a few pictures to convince herself.
\begin{lem}
Let $\alpha \in PSL_2\R$ be an isometry and $q_1 \neq q_2$ be two points in $\hyp^2$, neither of which is fixed by $\alpha$. Let $l$ denote the geodesic from $q_1$ to $q_2$, $l_1$ the geodesic segment from $q_1$ to $\alpha q_1$ and $l_2$ the geodesic segment from $q_2$ to $\alpha q_2$. Suppose we have two vector fields $\V_1, \V_2$ along $l_1,l_2$ respectively, and $c \in \Z$ such that:
\begin{enumerate}
\item $\V_1$ begins at $q_1$ pointing along $l$ towards $q_2$; ends at $\alpha q_1$ pointing along $\alpha l$ towards $\alpha q_2$; crosses the direction of $l_1$ transversely at finitely many points; the crossings taken with sign sum to $c$;
\item $\V_2$ begins at $q_2$ pointing along $l$ towards $q_1$; ends at $\alpha q_2$ pointing along $\alpha l$ towards $\alpha q_1$; crosses the direction of $l_2$ transversely at finitely many points; the crossings taken with sign sum to $c$.
\end{enumerate}
Then the vector fields $\V_1, \V_2$ both represent the same lift of $\alpha$ to $\widetilde{PSL_2\R}$.
\qed
\end{lem}

Since $\beta$ takes the segment $p_0 \rightarrow p_3$ to $p_1 \rightarrow p_2$, we may choose $\tilde{\beta}$ to be represented by a path of unit tangent vectors along the geodesic $p_0 \rightarrow p_1$ which begins pointing along $p_0 \rightarrow p_3$ and ends pointing along $p_1 \rightarrow p_2$; we may also choose $\tilde{\beta}$ to be represented by a path of unit tangent vectors along $p_3 \rightarrow p_2$, which begins pointing along $p_3 \rightarrow p_0$ and ends pointing along $p_2 \rightarrow p_1$. Using the above lemma, we may ensure that these two paths of tangent vectors represent the same $\tilde{\beta}$. Similarly, since $\alpha$ takes $p_4 \rightarrow p_1$ to $p_3 \rightarrow p_2$ we may choose $\tilde{\alpha}$ to be represented by the path of unit tangent vectors along the geodesic $p_4 \rightarrow p_3$ which begins pointing along $p_4 \rightarrow p_1$, ends pointing along $p_3 \rightarrow p_2$; and also we may choose $\tilde{\alpha}$ to be represented by tangent vectors along $p_1 \rightarrow p_2$ which begins pointing along $p_1 \rightarrow p_4$ and ends pointing along $p_2 \rightarrow p_3$. Chasing these paths of vectors around the pentagon, then, we obtain $\Twist ([\tilde{\alpha}^{-1}, \tilde{\beta}^{-1}],p) = \Delta[\Pent(\alpha,\beta;p)]$ on the nose.
\end{Proof}
It follows from the above proof that, even if $\Pent(\alpha,\beta;p)$ does not bound an immersed disc, we may still follow unit vectors and obtain $\Twist([\alpha^{-1},\beta^{-1}],p) \equiv 3\pi - \sum \theta_i$ modulo $2\pi$, where $\theta_i$ are the various angles between segments of $\Pent(\alpha,\beta;p)$.

If $S$ is a punctured torus with a hyperbolic structure and totally geodesic boundary, then a pentagon $\Pent(\alpha, \beta; p)$ is a fundamental domain for $S$, where $\alpha,\beta$ are the holonomy of a meridian and longitude, for appropriate choice of $p$, and $\Pent(\alpha,\beta;p)$ bounds an embedded disc; thus the area $\Delta[S]$ of $S$ is the area of the pentagon. By Gauss-Bonnet this area is $2\pi$. Also $[\alpha, \beta]$ is the holonomy of the boundary curve. Choosing $p \in \hyp^2$ as the appropriate vertex of the fundamental domain, we obtain that $\Twist([\alpha^{-1}, \beta^{-1}]), p) = \pm 2\pi$, sign depending on orientation. This is as it should be, since the developing image of the boundary should be the axis of its holonomy; and we may conclude that $[\alpha, \beta] \in \Hyp_{\pm 1}$.

More generally, whenever we have a pentagon $\Pent(\alpha, \beta; p)$ which bounds an immersed disc, the pentagon extends to a developing map for a hyperbolic cone-manifold structure on $S$ with piecewise geodesic boundary and one corner point; and every punctured torus with a hyperbolic cone-manifold structure with no interior cone points and at most one corner point can be cut into a pentagon in this way. For a complete investigation of such hyperbolic cone-manifold structures on punctured tori and their holonomy representations, see \cite{Me10MScPaper1}.

Consider again the twist of a commutator; it may be expanded as a product, using composition lemma \ref{lem:composition}:
\[
\Twist \left( [\tilde{\alpha}^{-1}, \tilde{\beta}^{-1}], p \right) = \Twist \Big( \tilde{\beta}, p \Big) + \Twist \Big( \tilde{\alpha}, \beta p \Big) + \Twist \Big( \tilde{\beta}^{-1}, \alpha \beta p \Big) + \Twist \Big( \tilde{\alpha}^{-1}, \beta^{-1} \alpha \beta p \Big) - \Delta[p_0,p_1,p_2,p_3,p_4],
\]
where $\Delta[p_0, p_1, p_2, p_3, p_4]$ is the signed area of the pentagon formed by geodesic segments (giving oriented boundary) $p_0 \rightarrow p_1 \rightarrow p_2 \rightarrow p_3 \rightarrow p_4$. Clearly not both $[p_0, p_1, p_2, p_3, p_4]$ and $[p_0, p_4, p_1, p_2, p_3] = \Pent(\alpha,\beta;p)$ can bound embedded discs! If area is however understood as a Euclidean angle defect given by angles between succeeding segments, or by cutting into triangles and summing their signed area, the above is true in all cases. Setting $q = \alpha \beta p$ the four twists involved can be simplified:
\begin{align*}
\Twist \Big( \tilde{\beta}, \beta^{-1} \alpha^{-1} q \Big) &+ \Twist \Big( \tilde{\alpha}, \alpha^{-1} q \Big) + \Twist \Big( \tilde{\beta}^{-1}, q \Big) + \Twist \Big( \tilde{\alpha}^{-1}, \beta^{-1} q \Big) \\
&= \Twist \Big( \tilde{\beta}, \alpha^{-1} q \Big) + \Twist \Big( \tilde{\alpha}, q \Big) - \Twist \Big( \tilde{\beta}, q \Big) - \Twist \Big( \tilde{\alpha}, \beta^{-1} q \Big).
\end{align*}
Here we have used equation \eqref{eq:translation} and lemma \ref{lem:inverse}. This can be considered as some kind of ``twist cross ratio'', the change in the twist of $\tilde{\alpha}$ under translation by $\beta$, relative to the change in twist of $\tilde{\beta}$ under translation by $\alpha$. With proposition \ref{prop:pentagon_twist}, these remarks give
\[
\Big( \Twist \big( \tilde{\alpha}, q \big) - \Twist \big( \tilde{\alpha}, \beta^{-1} q \big) \Big) - \Big( \Twist \big( \tilde{\beta}, q \big) - \Twist \big( \tilde{\beta}, \alpha^{-1} q \big) \Big) = \Delta[p_0, p_1, p_2, p_3, p_4] + \Delta[p_0, p_4, p_1, p_2, p_3].
\]

\subsection{Commutators and twist bounds}

Since a hyperbolic pentagon has area $< 3 \pi$, proposition \ref{prop:pentagon_twist} implies that $\left| \Twist \left( [\tilde{\alpha}, \tilde{\beta}],p \right) \right| < 3\pi$ when $\Pent(\alpha, \beta; p)$ is sufficiently nice. Such an inequality is true in general, as we now see.

Using the addition lemma \ref{lem:addition} gives two expressions for $\Twist \left( \tilde{\alpha}, p \right) + \Twist \left( \tilde{\beta}, p \right)$:
\[
\Twist \left( \tilde{\alpha}, p \right) + \Twist \left( \tilde{\beta}, p \right) = \Twist \left( \tilde{\beta} \tilde{\alpha}, p \right) + \Delta[p, \beta p, \beta \alpha p] = \Twist \left( \tilde{\alpha} \tilde{\beta}, p \right) + \Delta[p, \alpha p, \alpha \beta p].
\]
It follows that
\[
\Twist \left( \tilde{\alpha} \tilde{\beta}, p \right)  - \Twist \left( \tilde{\beta} \tilde{\alpha}, p \right) = \Delta[p, \beta p, \beta \alpha p] - \Delta[p, \alpha p, \alpha \beta p] .
\]
Applying the addition and inverse lemmas to the commutator $[\tilde{\alpha}, \tilde{\beta}] = \left( \tilde{\alpha} \tilde{\beta} \right) \left( \tilde{\beta} \tilde{\alpha} \right)^{-1}$, we have
\[
\Twist \left( [\tilde{\alpha}, \tilde{\beta}], p \right) = \Twist \left( \tilde{\alpha} \tilde{\beta}, p \right) - \Twist \left( \tilde{\beta} \tilde{\alpha}, p \right) - \Delta \left[ p, \alpha \beta p, [\alpha, \beta] p \right]
\]
Putting these together, we immediately have the following.
\begin{lem}[Commutator area]
\label{lem:commutator_area_formula}
For any $\tilde{\alpha}, \tilde{\beta} \in \widetilde{PSL_2\R}$ and $p \in \overline{\hyp}^2$,
\[
\Twist \left( \left[ \tilde{\alpha}, \tilde{\beta} \right], p \right) = \Delta[p, \beta p, \beta \alpha p] - \Delta[p, \alpha p, \alpha \beta p] - \Delta \left[ p, \alpha \beta p, [\alpha, \beta] p \right]
\]
\qed
\end{lem}
Hence for any $\tilde{\alpha}, \tilde{\beta} \in \widetilde{PSL_2\R}$ and $p \in \hyp^2$, we have (for $p \in S_\infty^1$ the inequality is not strict):
\begin{equation}
\left| \Twist \left( \left[ \tilde{\alpha}, \tilde{\beta} \right], p \right) \right| < 3\pi.
\label{eq:commutator_bounds}
\end{equation}

We can say more about the possible values for commutators; we consider the elliptic, parabolic, identity and hyperbolic cases separately.

If $[\alpha, \beta]$ is elliptic then take $p = \Fix [\alpha, \beta]$. The triangle formed by $p, \alpha \beta p, [\alpha, \beta] p$ then has zero area, and the twist of $[\tilde{\alpha}, \tilde{\beta}]$ at this point is $<2\pi$ in magnitude. Thus $[\tilde{\alpha}, \tilde{\beta}] \in \Ell_{-1} \cup \Ell_1$.

If $[\alpha, \beta]$ is parabolic, set $p = \Fix [\alpha, \beta] \in S_\infty^1$; then again $\Delta \left[ p, \alpha \beta p, [\alpha \beta] p  \right] = 0$, and so $\left| \Twist \left( \left[ \alpha, \beta \right], p \right) \right| \leq 2\pi$. In fact this twist must be in $\{-2\pi, 0, 2\pi\}$, and $[\tilde{\alpha}, \tilde{\beta}] \in \Par_{-1} \cup \Par_0 \cup \Par_1$ respectively.

\begin{figure}[tbh]
\centering
\includegraphics[scale=0.4]{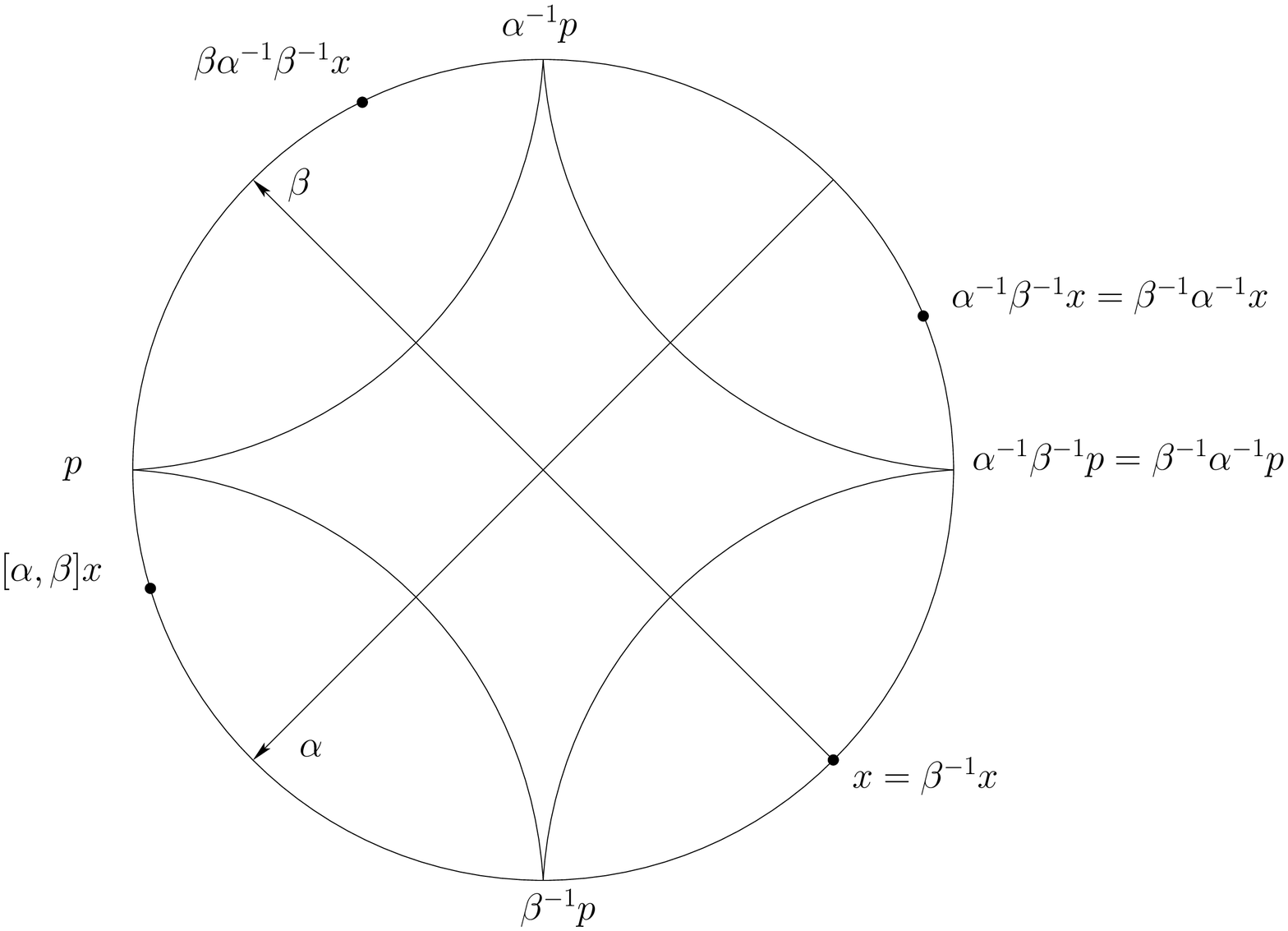}
\caption{Arrangement of axes in the case $[\alpha, \beta] \in \Par_1$.}
\label{fig:c}
\end{figure}

We can say something more in this case, with a little more work. Suppose that $[\tilde{\alpha}, \tilde{\beta}] \in \Par_1$; we will show in fact it lies in $\Par_1^-$. We have $\Twist \left( \left[ \alpha, \beta \right], p \right) = \Delta[ p, \beta p, \beta \alpha p ] - \Delta [ p, \alpha p, \alpha \beta p ]$. Applying the isometries $\alpha^{-1} \beta^{-1}$ and $\beta^{-1} \alpha^{-1}$ respectively to these ideal triangles gives this twist as $\Delta[\alpha^{-1} \beta^{-1} p, \alpha^{-1} p, p] - \Delta [ \beta^{-1} \alpha^{-1} p, \beta^{-1} p, p]$; these are both ideal triangles, and the values are $\pm \pi$ according to orientation, or $0$ if degenerate. Note that since $p = \Fix [\alpha, \beta]$, $\alpha^{-1} \beta^{-1} p = \beta^{-1} \alpha^{-1} p$. The only way to obtain $2\pi$ for the twist (and hence to lie in $\Par_1$), then, is if the four points $p, \beta^{-1} p, \alpha^{-1} \beta^{-1} p = \beta^{-1} \alpha^{-1} p, \alpha^{-1} p$ lie in anticlockwise order around the circle. Since $\alpha$ takes $(\alpha^{-1} p, \alpha^{-1} \beta^{-1} p) \mapsto (p, \beta^{-1} p)$, in different directions around $S_\infty^1$, $\alpha$ must be hyperbolic; similarly for $\beta$, and their axes must cross, and lie as shown in figure \ref{fig:c}. Let $x$ be the endpoint of $\Axis \beta$ shown; we now chase $x$ around the diagram to $[\alpha, \beta] x$. First $\beta^{-1} x = x$. As $\beta^{-1} x$ is anticlockwise of $\beta^{-1} p$, then $\alpha \beta^{-1} x$ must lie on the same side of $\Axis \alpha$ and anticlockwise of $\alpha^{-1} \beta^{-1} p$, hence as shown. Then $\alpha^{-1} \beta^{-1} x$ is anticlockwise of $\beta^{-1} \alpha^{-1} p$, so $\beta \alpha^{-1} \beta^{-1} x$ must lie anticlockwise of $\alpha^{-1} p$ and on the same side of $\Axis \beta$, hence as shown. By similar reasoning $[\alpha, \beta]x$ lies on the same side of $\Axis \alpha$ as $\beta \alpha^{-1} \beta^{-1} x$ and anticlockwise of $p$. Considering the result of applying $[\alpha, \beta]$ to $x$, we conclude that $[\alpha, \beta] \in \Par_1^{-}$. 

If $[\tilde{\alpha}, \tilde{\beta}] \in \Par_{-1}$, then the same argument applies and in fact it lies in $\Par_{-1}^+$.

If $[\alpha, \beta]$ is the identity, then we see that $[\tilde{\alpha}, \tilde{\beta}]$ is the identity in $\widetilde{PSL_2\R}$. If either of $\alpha$ or $\beta$ is the identity this is immediate; if not, $\alpha$ or $\beta$ are of the same type with same fixed points; and hence taking lifts and following unit vectors it is clear.

If $[\alpha, \beta]$ is hyperbolic, then equation \eqref{eq:commutator_bounds} immediately gives $[\tilde{\alpha}, \tilde{\beta}] \in \Hyp_{-1} \cup \Hyp_0 \cup \Hyp_1$. As it turns out, these are all possible.

\begin{figure}[tbh]
\centering
\includegraphics[scale=0.6]{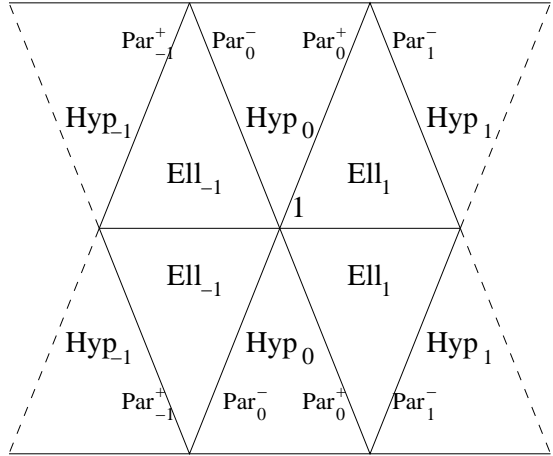}
\caption{Possible commutators in $\widetilde{PSL_2\R}$.}
\label{fig:7a}
\end{figure}

We have now proved the following theorem, which appears in \cite{Wood, Eisenbud-Hirsch-Neumann, Goldman88}; we also give a different proof in \cite{Mathews05}. See also figure \ref{fig:7a}.
\begin{thm}
\label{commutator_regions}
For $\tilde{\alpha}, \tilde{\beta} \in \widetilde{PSL_2\R}$,
\[
\left[ \tilde{\alpha}, \tilde{\beta} \right] \in \{1\} \cup \left( \bigcup_{n=-1}^1 \Hyp_n \cup \Ell_n \right) \cup \Par_0 \cup \Par_{-1}^+ \cup \Par_1^-.
\]
\qed
\end{thm}
Here we take $\Ell_0 = \emptyset$ for convenience. Combining this proposition with the trace lemma \ref{trace_regions} gives an immediate corollary.

\begin{cor}
\label{commutator_trace_info}
    If $\alpha,\beta \in PSL_2\R$ then
    \begin{enumerate}
        \item
            $\Tr[\alpha,\beta] > 2$ implies $[\alpha,\beta] \in \Hyp_0$;
        \item
            $\Tr[\alpha,\beta] = 2$ implies $[\alpha,\beta] \in \{1\} \cup \Par_0$;
        \item
            $\Tr[\alpha,\beta] \in (-2,2)$ implies $[\alpha,\beta] \in \Ell_{-1} \cup \Ell_1$;
        \item
            $\Tr[\alpha,\beta] = -2$ implies $[\alpha,\beta] \in \Par_{-1}^+ \cup \Par_1^-$;
        \item
            $\Tr[\alpha,\beta] < -2$ implies $[\alpha,\beta] \in \Hyp_{-1} \cup \Hyp_1$.
    \end{enumerate}\
\qed
\end{cor}

\subsection{Commutators and arrangements of axes}

In the proof of proposition \ref{commutator_regions} we showed that if $[\tilde{\alpha}, \tilde{\beta}] \in \Par_{\pm 1}^{\mp}$ then $\alpha, \beta$ are hyperbolic and their axes cross. We continue such analysis for other commutators, and make conclusions about the type and location of $\tilde{\alpha}, \tilde{\beta}$. In particular we prove the following result, which appears in \cite{Goldman03}.

\begin{prop}[Goldman \cite{Goldman03}]
\label{prop:axes_crossing}
Let $\alpha, \beta \in PSL_2\R$. The following are equivalent:
\begin{enumerate}
\item 
$\alpha, \beta$ are hyperbolic and their axes cross;
\item
$[\alpha, \beta] \in \Ell_{\pm 1} \cup \Par_{\pm 1}^{\mp} \cup \Hyp_{\pm 1}$;
\item
$\Tr [\alpha, \beta] < 2$.
\end{enumerate}
\end{prop}

We use a couple of lemmas. The first was implicitly used in the argument of the previous section and is straightforward. The second is a simple computation, for instance using Fermi coordinates (see e.g. \cite{Buser} p. 38).
\begin{lem}
\label{lem:hyperbolic_axis_placement}
Let $A,B,C,D$ be points on $S_\infty^1$ in anticlockwise order. Suppose $\alpha \in PSL_2\R$ is an isometry which takes $A$ to $B$ and $D$ to $C$. Then $\alpha$ is hyperbolic; the repulsive fixed point of $\alpha$ lies in the interval of $S_\infty^1$ between $D$ and $A$; and the attractive fixed point between $B$ and $C$.
\qed
\end{lem}

\begin{lem}
\label{lem:hyperbolic_translation_distance}
Let $\alpha \in PSL_2\R$ be hyperbolic and $x_1, y_1 \in \hyp^2$. The translation distances of $\alpha$ at $x_1$ and $x_2$ are equal, i.e. $d(x_1, \alpha x_1) = d(x_2, \alpha x_2)$, if and only if $x_1,x_2$ lie at the same perpendicular distance from $\Axis \alpha$.
\qed
\end{lem}

\begin{Proof}[Of proposition \ref{prop:axes_crossing}]
The equivalence of (ii) and (iii) is immediate from corollary \ref{commutator_trace_info}. To prove (ii) implies (i), we consider the various possible cases for $[\tilde{\alpha}, \tilde{\beta}]$.
\begin{itemize}
\item $[\tilde{\alpha}, \tilde{\beta}] \in \Hyp_{\pm 1}$.
The argument is virtually identical to the $\Par_{\pm 1}^\mp$ case. Consider $[\tilde{\alpha}, \tilde{\beta}] \in \Hyp_1$; the case $\Hyp_{-1}$ is identical with reversed orientation. Apply lemma \ref{lem:commutator_area_formula} taking $p \in S_\infty^1$ to be a fixed point of $[\tilde{\alpha}, \tilde{\beta}]$; so $\alpha^{-1} \beta^{-1} p = \beta^{-1} \alpha^{-1} p$ and, since $[\tilde{\alpha}, \tilde{\beta}] \in \Hyp_1$, $\Twist ( [\tilde{\alpha}, \tilde{\beta}], p) = 2\pi$. Note that $\Delta[p, \beta p, \beta \alpha p] = \Delta[\alpha^{-1} \beta^{-1} p, \alpha^{-1} p, p]$ and $\Delta[p, \alpha p, \alpha \beta p] = \Delta[\beta^{-1} \alpha^{-1} p, \beta^{-1} p, p]$. Then $p, \beta^{-1} p, \alpha^{-1} \beta^{-1} p = \beta^{-1} \alpha^{-1} p, \alpha^{-1} p$ must occur in anticlockwise order around $S_\infty^1$. By two applications of lemma \ref{lem:hyperbolic_axis_placement} then $\alpha, \beta$ are hyperbolic and their axes cross.

\item $[\tilde{\alpha}, \tilde{\beta}] \in \Par_{\pm 1}^{\mp}$.
We considered this case above, and concluded $\alpha, \beta$ hyperbolic with axes crossing.

\item $[\tilde{\alpha}, \tilde{\beta}] \in \Ell_{\pm 1}$.
Consider $[\tilde{\alpha}, \tilde{\beta}] \in \Ell_1$; the $\Ell_{-1}$ case is identical with reversed orientation. Apply lemma \ref{lem:commutator_area_formula} as above, taking $p = \Fix [\alpha, \beta]$. So $\alpha^{-1} \beta^{-1} p = \beta^{-1} \alpha^{-1} p$ and we have $\Twist ([\tilde{\alpha}, \tilde{\beta}],p) = \Delta[\alpha^{-1} \beta^{-1} p, \alpha^{-1} p, p] - \Delta[\beta^{-1} \alpha^{-1} p, \beta^{-1} p, p]$. These two triangles are \emph{congruent}: they both contain the side $p \rightarrow \alpha^{-1} \beta^{-1} p$, the isometry $\alpha$ takes the side $\alpha^{-1} \beta^{-1} p \rightarrow \alpha^{-1} p$ to $\beta^{-1} p \rightarrow p$, and the isometry $\beta$ takes the side $\beta^{-1} \alpha^{-1} p \rightarrow \beta^{-1} p$ to $\alpha^{-1} p \rightarrow p$. As $\Twist([\tilde{\alpha}, \tilde{\beta}],p) > 0$, the two congruent triangles fit together to give a non-self-intersecting quadrilateral $Q$ formed with vertices, in anticlockwise order, $p, \beta^{-1} p, \alpha^{-1} \beta^{-1} p = \beta^{-1} \alpha^{-1} p, \alpha^{-1} p$. Being constructed out of two congruent triangles, the opposite interior angles of $Q$ are equal. Moreover, extending the four side segments out to infinity, the only intersection points are the four vertices of the quadrilateral.

Extending the opposite sides $p \rightarrow \beta^{-1} p$ and $\alpha^{-1} p \rightarrow \alpha^{-1} \beta^{-1} p$ to infinity, then, they do not intersect; and they are related by $\alpha$. By lemma \ref{lem:hyperbolic_axis_placement} then $\alpha$ is hyperbolic, and $\Axis \alpha$ has fixed points at infinity separated by these lines. In particular, $\Axis \alpha$ intersects $p \rightarrow \beta^{-1} p$ (possibly extended) at a single point, and intersects $\alpha^{-1} p \rightarrow \alpha^{-1} \beta^{-1} p$ (possibly extended) at a single point also.  Moreover, since $\alpha^{-1} p \rightarrow p$ and $\beta^{-1} \alpha^{-1} p \rightarrow \beta^{-1} p$ have the same length (being related by $\beta$), the distance between $x$ and $\alpha(x)$ is the same for $x = \alpha^{-1} p$ and $x = \alpha^{-1} \beta^{-1} p$. By lemma \ref{lem:hyperbolic_translation_distance}, $\alpha^{-1} p$ and $\alpha^{-1} \beta^{-1} p$ lie at the same perpendicular distance from $\Axis \alpha$. Hence they lie on opposite sides of $\Axis \alpha$, and $\Axis \alpha$ intersects the segment $\alpha^{-1} p \rightarrow \beta^{-1} \alpha^{-1} p$. By the same argument regarding translation distances of $\alpha^{-1}$, $\Axis \alpha^{-1} = \Axis \alpha$ intersects the segment $p  \rightarrow \beta^{-1} p$. That is, $\Axis \alpha$ intersects two opposite sides of the quadrilateral $Q$. By the same argument, $\beta$ is hyperbolic and $\Axis \beta$ intersects the other pair of opposite sides of $Q$. Hence $\Axis \alpha$ and $\Axis \beta$ intersect.
\end{itemize}

To prove (i) implies (iii), we repeat the argument of \cite{Goldman03}, writing matrices for $\alpha, \beta \in PSL_2\R$. We may conjugate so that $\alpha$ has fixed points $-1,1$ in the upper half plane, and $\beta$ has fixed points $r,\infty$ where $r \in (-1,1)$. We may write
\[
\alpha = \begin{pmatrix} \cosh x & \sinh x \\ \sinh x & \cosh x \end{pmatrix}, \quad
\beta = \begin{pmatrix} e^y & - 2r \sinh y \\ 0 & e^{-y} \end{pmatrix}
\]
where $x, y \in \R$. A calculation then gives
\[
\Tr [\alpha, \beta] = 2 + 4 \left( r^2 - 1 \right) \sinh^2 x \sinh^2 y < 2.
\]
\end{Proof}

Note Goldman's proof in \cite{Goldman03} of this proposition is entirely algebraic; the difference here is that we have proved one direction geometrically.

\subsection{Surface group representations and the Milnor--Wood inequality}

Now consider $\tilde{\alpha}_1, \tilde{\beta}_1, \ldots, \tilde{\alpha}_g, \tilde{\beta}_g, \tilde{\gamma}_1, \ldots, \tilde{\gamma}_n \in \widetilde{PSL_2\R}$ and consider 
\[
\left[ \tilde{\alpha}_1, \tilde{\beta}_1 \right] \; \left[ \tilde{\alpha}_2, \tilde{\beta}_2 \right] \; \cdots \; \left[ \tilde{\alpha}_g, \tilde{\beta}_g \right] \tilde{\gamma}_1 \tilde{\gamma_2} \cdots \tilde{\gamma}_n.
\]
The commutators, as we saw in lemma \ref{commutator_lift}, are independent of choice of lift of $\alpha_i, \beta_i$; for the $\tilde{\gamma}_n$, let us assume they are ``efficiently chosen'', i.e. have twist less than $\pi$ in magnitude at some point $p$. Such an expression is of course the relator in the standard presentation of the fundamental group of the surface $S$ of genus $g$ with $n$ boundary components. Such $\alpha_i, \beta_i, \gamma_i \in PSL_2\R$ arise as the holonomy of a hyperbolic structure on $S$, and in this case the product is $1 \in PSL_2\R$; it therefore lifts to some $\z^m \in \widetilde{PSL_2\R}$. We ask how large $m$ can be, i.e. how large the twist of the relator expression can be.

Repeated use of the addition lemma gives a bound, in the ``hyperbolic'' case $\chi(S) = 2-2g-n < 0$.
\begin{thm}[Milnor \cite{Milnor}]
\label{bounds for Euler class on general surface}
Let $p \in \hyp^2$; let $\tilde{\alpha}_1, \tilde{\beta}_1, \ldots, \tilde{\alpha}_g, \tilde{\beta}_g, \tilde{\gamma}_1, \ldots, \tilde{\gamma}_n \in \widetilde{PSL_2\R}$ with $2 - 2g - n < 0$ and $\left| \Twist \left( \tilde{\gamma}_i , p \right) \right| \leq \pi$; assume
\[
\left[ \alpha_1, \beta_1 \right] \; \left[ \alpha_2, \beta_2 \right] \; \cdots \; \left[ \alpha_g, \beta_g \right] \gamma_1 \gamma_2 \cdots \gamma_n = 1 \in PSL_2\R.
\]
Then
\[
\left[ \tilde{\alpha}_1, \tilde{\beta}_1 \right] \; \left[ \tilde{\alpha}_2, \tilde{\beta}_2 \right] \; \cdots \; \left[ \tilde{\alpha}_g, \tilde{\beta}_g \right] \tilde{\gamma}_1 \tilde{\gamma_2} \cdots \tilde{\gamma}_n = \z^m
\]
where $|m| \leq |\chi(S)| = 2g+n-2$.
\end{thm}

\begin{Proof}
There are $4g+n$ terms in the expression; we use the addition lemma $4g+n-1$ times. Note that whenever $\tilde{\xi} \tilde{\eta} = 1$ we have $\Delta[p, \xi p, \xi \eta p] = 0$, being a degenerate triangle. So the first time we use the addition lemma the triangle is degenerate, and there are at most $4g+n-2$ nondegenerate triangles. Once we have used the addition lemma $4g+n-1$ times, we have the difference between the twist of the relator and the twist of the individual terms as a sum of signed areas of $4g+n-2$ triangles; denote these signed areas $\Delta_1, \ldots, \Delta_{4g+n-2}$. The inverses from the commutators cancel, thanks to lemma \ref{lem:inverse}, leaving only the twists of the $\tilde{\gamma}_i$. Thus
\[
\Twist \left( \z^m, p \right) = \sum_{i=1}^n \Twist \left( \tilde{\gamma}_i, p \right) + \sum_{i=1}^{4g+n-2} \Delta_i,
\]
and since triangle areas are $<\pi$ and the twists of the $\tilde{\gamma}_i$ are by assumption $\leq \pi$,
\[
2\pi \; |m| = \Big| \Twist \left( \z^m, p \right) \Big| < n \pi + (4g+n-2) \pi = 2\pi \left( 2g + n - 1 \right).
\]
Since $|m| < 2g+n-1$ and $m$ is an integer, we are done.
\end{Proof}

Note that the above argument does not work for the ``non-hyperbolic'' case $2g+n \leq 2$. For instance, with $k=0, n=2$, setting $\tilde{c}_1 = \tilde{c}_2$ to be half turns of twist $\pi$, we have $c_1 c_2 = 1 \in PSL_2\R$ but $\tilde{c}_1 \tilde{c}_2 = \z$.

Suppose we have a representation $\rho: \pi_1(S) \To PSL_2\R$. After choosing lifts $\tilde{\gamma}_i$ of the images of the boundary components, we have a lift of the image of the relator, which is some $\z^m$. This $m$ is essentially the (relative) \emph{Euler class} $\E(\rho)$ of the representation: see e.g. \cite{Goldman88,  Me10MScPaper2}. More precisely $\E(\rho) \in H^2(S)$ and takes the fundamental class $[S]$ to $m$. Interpreting the above result this way we have:
\begin{thm}
\label{Milnor-Wood}
Let $\rho: \pi_1(S) \To PSL_2\R$ be a representation. Taking lifts $\tilde{\gamma}_i$ of the images of the boundary components with $|\Twist(\tilde{\gamma}_i,p)|<\pi$ at some point $p$, the (relative) Euler class $\E(\rho)$ takes the fundamental class $[S]$ to $m \in \Z$, where $|m| \leq |\chi(S)|$.
\qed
\end{thm}

The above two results, and similar formulations, are known generally as the \emph{Milnor-Wood inequality}. The first inequality was first proved by Milnor \cite{Milnor}, generalised by Wood \cite{Wood}, reproved by Goldman \cite{Goldman88}, and generalised further by Eisenbud--Hirsch--Neumann \cite{Eisenbud-Hirsch-Neumann}. The reformulation in terms of Euler class was first given, so far as we know, by Wood \cite{Wood}.

\subsection{Larger products and polygons}

Using our existing results --- composition lemma \ref{lem:composition}, addition lemma \ref{lem:addition}, commutator pentagon \ref{prop:pentagon_twist}, and commutator area \ref{lem:commutator_area_formula} --- together, we obtain results for more complicated expressions and figures.

For instance, using the composition lemma repeatedly on a product $\tilde{\gamma}_1, \tilde{\gamma}_2 \cdots \tilde{\gamma}_n$ gives
\[
\Twist \left( \prod_{i=1}^n \tilde{\gamma}_i, p \right) = \sum_{i=1}^n \Twist \left( \tilde{\gamma}_i, \left( \prod_{j=i+1}^n \gamma_j \right) p \right) - \sum_{i=2}^n \Delta \left[ \left( \prod_{j=i+1}^n \gamma_j \right) p, \left( \prod_{j=i}^n \gamma_j \right) p, \left( \prod_{j=1}^n \gamma_j \right) p \right].
\]
But these triangles share successive sides. Consider the polygon $\CC({\bf \gamma}; p)$ with $n+1$ sides
\[
p \rightarrow \gamma_1 \gamma_2 \cdots \gamma_n p \rightarrow \gamma_2 \cdots \gamma_n p \rightarrow \cdots \rightarrow \gamma_{n-1} \gamma_n p \rightarrow \gamma_n p \rightarrow p,
\]
oriented by the direction on the boundary above. If $\CC({\bf \gamma};p)$ is convex, then it can be cut into precisely the triangles occurring in the above sum; and so, as long as the polygon is simple (i.e. non-self-intersecting), the above expression gives the area.
\begin{lem}[Composition polygon]
\label{lem:composition_polygon}
Let $\tilde{\gamma}_1, \ldots, \tilde{\gamma}_n \in \widetilde{PSL_2\R}$ and $p \in \overline{\hyp}^2$. Suppose the polygon $\CC(\gamma_1, \ldots, \gamma_n; p)$ is simple. Then
\[
\Twist \left( \prod_{i=1}^n \tilde{\gamma}_i, p \right) = \sum_{i=1}^n \Twist \left( \tilde{\gamma}_i, \left( \prod_{j=i+1}^n \gamma_j \right) p \right) + \Delta \left[ \CC \left( \gamma_1, \ldots, \gamma_n; p \right) \right].
\]
\qed
\end{lem}
Note that if $\gamma_1 \cdots \gamma_n = 1 \in PSL_2\R$ then $\CC(\gamma_1, \ldots, \gamma_n; p)$ reduces to an $n$-gon and the result still holds.

\begin{figure}[tbh]
\centering
\includegraphics[scale=0.45]{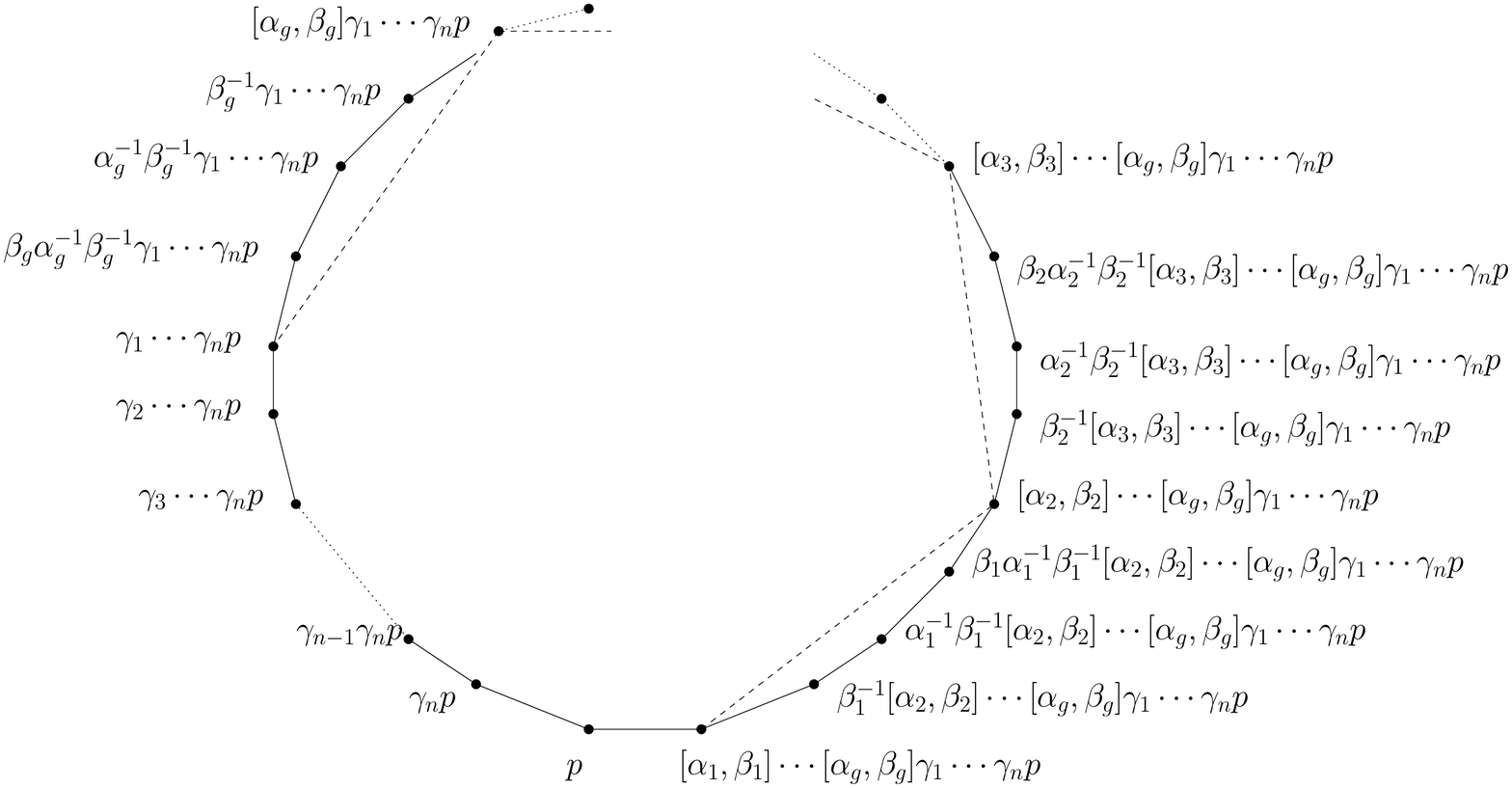}
\caption{The polygon $\D({\bf \alpha}, {\bf \beta}, {\bf \gamma};p)$.}
\label{fig:d}
\end{figure}

Consider now polygons of the type arising as fundamental domains of hyperbolic surfaces. Given $\tilde{\alpha_1}, \tilde{\beta}_1, \cdots, \tilde{\alpha}_g, \tilde{\beta}_g, \tilde{\gamma}_1, \ldots, \tilde{\gamma}_n$, we consider a polygon $\D({\bf \alpha}, {\bf \beta}, {\bf \gamma}; p)$ associated to the surface group relator
\[
\left[ \tilde{\alpha}_1, \tilde{\beta}_1 \right] \; \left[ \tilde{\alpha}_2, \tilde{\beta}_2 \right] \cdots \left[ \tilde{\alpha}_g, \tilde{\beta}_g \right] \tilde{\gamma}_1 \tilde{\gamma}_2 \cdots \tilde{\gamma}_n.
\]
The vertices of $\D({\bf \alpha}, {\bf \beta}, {\bf \gamma}; p)$ are, along the oriented boundary:
\begin{align*}
p &\rightarrow [\alpha_1, \beta_1] \cdots [\alpha_g, \beta_g] \gamma_1 \cdots \gamma_n p \\
& \rightarrow \beta_1^{-1} [\alpha_2, \beta_2] \cdots [\alpha_g, \beta_g] \gamma_1 \cdots \gamma_n p 
\rightarrow \alpha_1^{-1} \beta_1^{-1} [\alpha_2, \beta_2] \cdots \gamma_n p \rightarrow \beta_1 \alpha_1^{-1} \beta_1^{-1} [\alpha_2, \beta_2] \cdots \gamma_n p \\
& \rightarrow [\alpha_2, \beta_2] \cdots [\alpha_g, \beta_g] \gamma_1 \cdots \gamma_n p \\
& \rightarrow \beta_2^{-1} [\alpha_3, \beta_3] \cdots [\alpha_g, \beta_g] \gamma_1 \cdots \gamma_n p 
\rightarrow \alpha_2^{-1} \beta_2^{-1} [\alpha_3, \beta_3] \cdots \gamma_n p \rightarrow \beta_2 \alpha_2^{-1} \beta_2^{-1} [\alpha_3, \beta_3] \cdots \gamma_n p \\
& \rightarrow [\alpha_3 \beta_3] \cdots [\alpha_g, \beta_g] \gamma_1 \cdots \gamma_n p \\
& \rightarrow \cdots \\
& \rightarrow \gamma_1 \cdots \gamma_n p \\
& \rightarrow \gamma_2 \cdots \gamma_n p \rightarrow \gamma_3 \cdots \gamma_n p \rightarrow \cdots \rightarrow \gamma_n p \\
& \rightarrow p
\end{align*}
See figure \ref{fig:d}.

Note that $\D \left( {\bf \alpha}, {\bf \beta}, {\bf \gamma}; p \right)$ can be considered as a sequence of pentagons attached to the polygon 
\[
\CC([\alpha_1, \beta_1], [\alpha_2, \beta_2], \ldots, [\alpha_g, \beta_g], \gamma_1, \ldots, \gamma_n; p).
\]
(See the dashed lines in figure \ref{fig:d}.) Note that $\D({\bf \alpha}, {\bf \beta}, {\bf \gamma};p)$ is a $(4g+n+1)$-gon, but in a surface group representation $[\alpha_1, \beta_1] \cdots [\alpha_g, \beta_g] \gamma_1 \cdots \gamma_n = 1$, so it is a $(4g+n)$-gon. In any case, applying the composition polygon lemma \ref{lem:composition_polygon} gives
\begin{align*}
\Twist \left( \prod_{i=1}^g \left[ \tilde{\alpha}_i, \tilde{\beta}_i \right] \prod_{i=1}^n \tilde{\gamma}_i, p \right)
&=
\sum_{i=1}^g \Twist \left( \left[ \tilde{\alpha}_i, \tilde{\beta}_i \right], \left( \prod_{j=i+1}^g \left[ \alpha_j, \beta_j \right] \right) \left( \prod_{j=1}^n \gamma_j \right) p \right) \\
& \quad + \sum_{i=1}^n \Twist \left( \tilde{\gamma}_i, \left( \prod_{j=i+1}^n \tilde{\gamma}_j \right) p \right) + \Delta \left[ \CC \left( \left[ \alpha_1, \beta_1 \right], \ldots, \left[ \alpha_g, \beta_g \right], \gamma_1, \ldots, \gamma_n; p \right) \right].
\end{align*}
But the twist of a commutator is the area of a pentagon (proposition \ref{prop:pentagon_twist}), hence the above is equal to
\begin{align*}
\sum_{i=1}^n \Twist \left( \tilde{\gamma}_i, \left( \prod_{j=i+1}^n \tilde{\gamma}_j \right) p \right) + \Delta \left[ \CC \left( \left[ \alpha_1, \beta_1 \right], \ldots, \left[ \alpha_g, \beta_g \right], \gamma_1, \ldots, \gamma_n; p \right) \right] \\
+ 
\sum_{i=1}^g \Delta \left[ \Pent \left( \alpha_i^{-1}, \beta_i^{-1}; \left( \prod_{j=i+1}^g \left[ \alpha_j, \beta_j \right] \right) \left( \prod_{j=1}^n \gamma_j \right) p \right) \right]
\end{align*}
These areas simply add up to our polygon $\D(\alpha, \beta, \gamma; p)$ --- provided that we have simple polygons. In order to ensure that the entire polygon, and the result of cutting off these pentagons, are simple, the easiest thing to do is make an assumption of \emph{convexity}.
\begin{thm}
Let $\tilde{\alpha_1}, \tilde{\beta}_1, \cdots, \tilde{\alpha}_g, \tilde{\beta}_g, \tilde{\gamma}_1, \ldots, \tilde{\gamma}_n \in \widetilde{PSL_2\R}$ and $p \in \overline{\hyp}^2$. Suppose that $\D({\bf \alpha}, {\bf \beta}, {\bf \gamma}; p)$ is a convex polygon. Then
\[
\Twist \left( \left[ \tilde{\alpha}_1, \tilde{\beta}_1 \right] \cdots \left[ \tilde{\alpha}_g, \tilde{\beta}_g \right] \tilde{\gamma}_1 \cdots \tilde{\gamma}_n, p \right) = \sum_{i=1}^n \Twist \left( \tilde{\gamma}_i, \gamma_{i+1} \gamma_{i+2} \cdots \gamma_n  p \right) + \Delta \left[ \D \left( {\bf \alpha}, {\bf \beta}, {\bf \gamma}; p \right) \right].
\]
\qed
\end{thm}

This gives another proof of the Milnor--Wood inequality, provided $\D({\bf \alpha}, {\bf \beta}, {\bf \gamma}; p)$ is convex: if $\rho: \pi_1 (S) \To PSL_2\R$ is a representation and $\E(\rho)=m$ then
\[
2\pi m = \sum_{i=1}^n \Twist \left( \tilde{\gamma}_i, \gamma_{i+1} \gamma_{i+2} \cdots \gamma_n  p \right) + \Delta \left[ \D \left( {\bf \alpha}, {\bf \beta}, {\bf \gamma}; p \right) \right]
\]
As before, we choose lifts $\tilde{\gamma}_i$ of the boundary components with twists less than $\pi$ in magnitude. Noting that the area of a hyperbolic $(4g+n)$-gon is less than $(4g+n-2)\pi$, we immediately have
\[
2\pi |m| < n\pi + \pi \left( 4g+n-2 \right) = 2 \pi \left( 2g+n-1 \right),
\]
and again as $m$ is an integer, $|m| \leq 2g+n-2 = |\chi(S)|$.

\section{Milnor's angle function}
\label{sec:Milnors_function}

\subsection{Definition}

We now recall the definition of Milnor's angle function, relate it to our notion of twist, and then deduce various properties from our results on twisting.

A matrix $\alpha \in GL_2^+ \R$ can be written uniquely in the form $\alpha = R(\alpha) S(\alpha)$, where $R(\alpha),S(\alpha) \in GL_2^+\R$, $R$ is orthogonal (i.e. $R \in SO_2\R$) and $S$ is symmetric positive definite. Since $R$ is orthogonal, it is of the form
\[
    R(\alpha) = \begin{pmatrix} \cos \theta & \sin \theta \\ -
    \sin \theta & \cos \theta \end{pmatrix}
\]
for some $\theta$. This $\theta$ can be thought of as the angle of rotation of $\alpha$; and hence we have a function $\theta: GL_2^+ \R \To \R / 2 \pi \Z$. Milnor's $\Theta$ function is a lift of this map, $\Theta: \widetilde{GL_2^+ \R} \To \R$. 

Indeed, the map $R: SL_2\R \To SO_2\R \cong S^1$ is a retraction, which lifts to a retraction $\tilde{R}: \widetilde{SL_2\R} \To \widetilde{SO_2\R} \cong \R$. Since $SO_2\R \cong S^1$, we have $\widetilde{SO_2\R} \cong \R$, $\mathfrak{so}_2\R \cong \R$, so the
exponential map is
\[
    \exp: \R \To SO_2\R, \qquad \theta \mapsto \begin{pmatrix}\cos
    \theta & \sin \theta \\ - \sin \theta & \cos \theta
    \end{pmatrix}.
\]
lifts to $\widetilde{\exp}: \R \To \widetilde{SO_2\R} \subset
\widetilde{PSL_2\R}$. We define the angle
function $\Theta: \widetilde{PSL_2\R} \To \R$ by
\[
    \Theta(\tilde{\alpha}) = \widetilde{\exp}^{-1} \left(
    \tilde{R}(\alpha)\right).
\]

Although Milnor considers $\Theta$ on $\widetilde{GL_2^+ \R}$ and general $2 \times 2$ real matrices with positive determinant, we only consider $\widetilde{SL_2\R}$ and $2 \times 2$ real matrices with determinant $1$. But any matrix with positive determinant can act as a hyperbolic isometry on the upper half plane as a fractional linear transformation; indeed, as groups, $GL_2^+ \R \cong \R_+ \times SL_2\R$, under the isomorphism $A \mapsto (\det A, A / \det A)$, and the fractional linear transformations of $A$ and $\lambda A$ are equal for any $\lambda \in \R_+$. We see that $\Theta(\lambda A) = \Theta (A)$. So restricting to $SL_2\R$ in essence loses no generality; if we like we could define twist on $\widetilde{GL_2^+ \R} \times \overline{\hyp}^2$ without any difficulty.

\subsection{A geometric interpretation}

We now interpret Milnor's $\Theta$ through twisting. Although it seems that similar ideas have been used previously, for instance in \cite{Goldman88, Wood}, as far as we know this has not been described explicitly before. 
\begin{prop}
\label{geometric_interpretation}
    Let $\tilde{\alpha} \in \widetilde{PSL_2\R}$. Then
    \[
        \Theta(\tilde{\alpha}) = \frac{1}{2}
        \Twist \left( \tilde{\alpha}, i \right).
    \]
\end{prop}

For the proof, we begin with some simple observations.
\begin{lem}
\label{lem:passes_through_i}
    In the upper half plane model, the geodesic with endpoints at
    infinity $a,b \in \R$ passes through $i$ if and only if $ab=-1$.
\qed
\end{lem}

\begin{lem}
\label{rotation_about_i}
In the upper half plane model, the matrix
\[
R(\theta) = \begin{pmatrix} \cos \theta & - \sin \theta \\ \sin \theta & \cos \theta \end{pmatrix} \in SO_2\R \subset SL_2\R,
\]
taken modulo sign (i.e. in $PSL_2\R = \Isom^+ \hyp^2$), acts as a rotation of angle $2\theta$ anticlockwise about $i$.
\qed
\end{lem}

\begin{lem}
\label{symmetric positive definite result}
    An isometry of $\hyp^2$ is represented by
    a symmetric positive definite matrix in $SL_2\R$ other
    than the identity if and only
    if it is hyperbolic and its axis passes through $i$.
\end{lem}

\begin{Proof}
A simple computation, contained in \cite{Mathews05}.
\end{Proof}

\begin{Proof}[of proposition \ref{geometric_interpretation}]
    From the definition of $\Theta$ we have
    \[
        \tilde{\alpha} =
        \widetilde{\exp} \Big( \Theta(\tilde{\alpha}) \Big)
        \tilde{S}(\tilde{\alpha})
    \]
    where $\widetilde{\exp} \left( \Theta ( \tilde{\alpha} ) \right) = \tilde{R}(\tilde{\alpha})$ and $\tilde{S} (\tilde{\alpha})$ is the lift of a symmetric positive definite matrix; since $\tilde{R}$ is a retraction, $\tilde{S}(\tilde{\alpha})$ is connected to the identity through lifts of symmetric positive definite matrices (hence lifts of hyperbolic isometries). Thus $\tilde{S}(\tilde{\alpha}) \in \Hyp_0 \cup \{1\}$.

    Consider now the action of $\tilde{\alpha}$ on the hyperbolic plane. Since $S(\tilde{\alpha})$ is symmetric positive definite, by \ref{symmetric positive definite result} it is either the identity, or a translation along an axis passing through $i$. And $\tilde{S}(\tilde{\alpha}) \in \Hyp_0 \cup \{1\}$ is the simplest lift of this translation. This action is followed by that of $\widetilde{\exp}\left( \Theta(\tilde{\alpha}) \right)$, which is a rotation of angle $2 \Theta(\tilde{\alpha})$ about $i$ (by \ref{rotation_about_i}). So the overall action of $\tilde{\alpha}$ is to translate from $\alpha^{-1}(i)$ to $i$, and then rotate by angle $2 \Theta(\tilde{\alpha})$. Thus $\Twist \left( \tilde{\alpha}, \alpha^{-1}(i) \right) = 2 \Theta(\tilde{\alpha})$. By equation \ref{eq:translation}, this is also equal to $\Twist \left( \tilde{\alpha}, i \right)$.
\end{Proof}

We note a comparison with Wood's methods in \cite{Wood}, which presage subsequent work on the Milnor--Wood inequality for homeomorphisms of the circle (e.g.  \cite{Eisenbud-Hirsch-Neumann, Calegari04}). As mentioned above in section \ref{sec:deriv_isoms}, Wood regards elements of $SL_2\R$ as acting on the $S^1$ of oriented lines through the origin in $\R^2$, which is equivalent to considering the action on the circle at infinity; we consider the action on the $S^1$ of unit tangent vectors at a point. Wood shows that his function $\widetilde{\Top^+ S^1} \To \R$ is an extension of Milnor's $\Theta$. We have shown that our twist at the particular point $i$ is equal to $\Theta$. Hence the action we consider on unit tangent vectors at $i$ is equivalent to Wood's action on the circle at infinity. Although our function $\Twist$ is only defined on group $\widetilde{PSL_2\R}$, which is much smaller than $\widetilde{\Top^+ S^1}$, in applying to different points in $\hyp^2$ it is in a sense a generalisation of Wood's function.

\subsection{Properties}

We can now immediately deduce properties of $\Theta$ from $\Twist$. Note several of our results for $\Twist$, particularly involving surface group relators, cannot be expressed simply in terms of $\Theta$, since they involve twists at different points. In this sense, the twist is a true generalisation of $\Theta$.

Immediately from proposition \ref{twist_bounds}, we have the values of $\Theta$ on the various regions of $\widetilde{PSL_2\R}$:
\begin{align*}
        \Theta(\Hyp_n) &= \left( \left(n-\frac{1}{2} \right) \pi,
        \left( n + \frac{1}{2} \right) \pi \right) \\
        \Theta(\Par_n) &= \left( \left(n-\frac{1}{2} \right) \pi,
        \left( n + \frac{1}{2} \right) \pi \right) \\
        \Theta(\Ell_n) &= \left\{ \begin{array}{ll} \Bigl( (n-1)\pi, n\pi \Bigr) & \text{
        for $n>0$} \\
        \Bigl( -|n|\pi, (-|n|+1)\pi \Bigr) &
        \text{ for $n<0$} \end{array} \right.
\end{align*}\

From the inverse lemma \ref{lem:inverse}, we immediately obtain $\Theta(\alpha) = - \Theta (\alpha^{-1})$. This also follows from the definition of $\Theta$. We have $R(\alpha^{-1}) = R(\alpha)^{-1}$, and by continuity we can deduce $\tilde{R}(\tilde{\alpha}^{-1}) = \tilde{R}(\tilde{\alpha})^{-1}$, so that $\Theta(\tilde{\alpha}^{-1}) = - \Theta(\tilde{\alpha})$.

From proposition \ref{prop:pentagon_twist} we have that if $\Pent(\alpha,\beta;i)$ bounds an embedded (or immersed) disc, then
\begin{equation}
\label{eq:theta_pentagon}
\Theta \left( \left[ \tilde{\alpha}^{-1}, \tilde{\beta}^{-1} \right] \right) = \frac{1}{2} \Delta \left[ \Pent \left( \alpha, \beta; i \right) \right].
\end{equation}
From lemma \ref{lem:commutator_area_formula} we have:
\begin{equation}
2 \Theta \left( \left[ \tilde{\alpha}, \tilde{\beta} \right] \right) = \Delta[i, \beta i, \beta \alpha i] - \Delta[i, \alpha i, \alpha \beta i] - \Delta \left[ i, \alpha \beta i, [\alpha, \beta] i \right], \quad \text{hence} \quad 
\left| \Theta \left( \left[ \tilde{\alpha}, \tilde{\beta} \right] \right) \right| < \frac{3\pi}{2}.
\end{equation}
The quasimorphism property of $\Theta$ follows immediately from the addition lemma \ref{lem:addition}.
\begin{equation}
\label{eq:theta_quasimorphism}
\Theta \left( \tilde{\beta} \tilde{\alpha} \right) = \Theta \left( \tilde{\beta} \right) + \Theta \left( \tilde{\alpha} \right) - \frac{1}{2} \Delta[i, \beta i, \beta \alpha i], \quad \text{hence} \quad 
\left| \Theta \left( \tilde{\beta} \tilde{\alpha} \right) - \Theta \left( \tilde{\beta} \right) - \Theta \left( \tilde{\alpha} \right) \right| < \frac{\pi}{2}
\end{equation}

Repeated application of the final inequality is enough to prove the Milnor--Wood inequality, as Milnor carries out in \cite{Milnor}. But the equalities involving areas in \eqref{eq:theta_pentagon}--\eqref{eq:theta_quasimorphism} appear to be new.

\addcontentsline{toc}{section}{References}

\small

\bibliography{danbib}
\bibliographystyle{amsplain}

\end{document}